\documentclass{article}
\usepackage[english, francais]{babel}
\usepackage{amsfonts}
\usepackage[T1]{fontenc}
\usepackage{graphicx}
\usepackage{color}
\newcommand{\pre}{{\noindent\noindent{} \bf Preuve. }}
\newcommand{\di}{\displaystyle}
\def\sk{\vskip 0.2cm}
\def\Z{\mathbf Z}
\def\N{{\mathbf N}}
\def\squareforqed{\hbox{\rlap{$\sqcap$}$\sqcup$}}
\def\qed{\ifmmode\else\unskip\quad\fi\squareforqed}
\def\smartqed{\def\qed{\ifmmode\squareforqed\else{\unskip\nobreak\hfil
\penalty50\hskip1em\null\nobreak\hfil\squareforqed
\parfillskip=0pt\finalhyphendemerits=0\endgraf}\fi}}
\def\smartqed{\def\qed{\ifmmode\squareforqed\else{\unskip\nobreak\hfil
\penalty50\hskip1em\null\nobreak\hfil\squareforqed
\parfillskip=0pt\finalhyphendemerits=0\endgraf}\fi}}
\newcommand{\cqfd}{\textcolor{blue}{\smartqed \qed}}
\newcommand{\al}{\alpha}

\newcommand{\be}{\beta}
\newcommand{\ep}{\varepsilon}

\newcommand{\Om}{\Omega}

\newcommand{\te}{\theta}

\newtheorem{thm}{Th\'eor\`eme}[section]
\newtheorem{cor}[thm]{Corollaire}
\newtheorem{lem}[thm]{Lemme}
\newtheorem{pro}[thm]{Proposition}

\newtheorem{exe}[thm]{Exemple}

\newtheorem{df}{D\'efinition}[section]
\newtheorem{rem}{Remarque}[section]

\title{Extensions $+\infty$-$w_0$-g\'en\'er\'ees}
\author{El Hassane Fliouet}
\date{}
\begin{document}
%\date{Received: date / Accepted: date}
\maketitle
\markboth{E. Fliouet}{Extensions $+\infty$-$w_0$-g\'en\'er\'ees}
%%%%%%%%%%%%%%%%%%%%%%%%%%%%%%%%%%%%%%%%%%%%%%%%%%%%%%%
%%%%%%%%%%%%%%%%%%%%%%%%%%%%%%%%%%%%%%%%%
\begin{abstract}
In this note, we continue to be interested in the relationship that connects  the restricted distribution of finitude at the local level of intermediate fields of a purely inseparable extension $K/k$ to the absolute or global finitude of $K/k$.  In "{\it $w_0$-generated field extensions,}
Arch. Math. {\bf  47}, (1986), 410-412",  JK Deveney constructed an example of modular extension $K/k$   called $w_0 $-generated   such that for any proper subfield $L$ of $K/k $, $L$ is finite over $k$, and  for every $ n \in \N$, we have $ [k^{p^{- n}} \cap K: k] = p^{2n} $. This example has proved to be extremely useful in the construction of other examples of $w_0$-generated extensions. In particular, we prolong the $w_0$-generated to an extension of unspecified finite size.
However, when $K/k$ is of unbounded size, we show that any modular extension of unbounded exponent admits a  proper subextension of unbounded exponent.  This brings us  to study the $w_0$-generated in the restricted sense. In addition, with the aim of extending the $w_0$-generated to a purely inseparable extension of unbounded size, we propose other generalizations.
\end{abstract}
{\bf Mathematics Subject Classification MSC2010:} Primary 12F15
\sk

\selectlanguage{francais}
\section{Introduction}
Soit $K/k$ une extension purement ins\'{e}parable de caract\'{e}ristique $p>0$. Une partie $B$ de $K$ est dite $r$-base de $K/k$ si $K=k(K^p)(B)$, et pour tout $x\in B$, $x\not\in k(K^p)(B\setminus \{x\})$. En vertu du (\cite{N.B}, III, p. 49, corollaire \textcolor{blue}{3}) et de la propri\'{e}t\'{e} d'\'{e}change des $r$-bases,  on en d\'{e}duit que toute extension admet une $r$-base et que le cardinal d'une $r$-base  est invariant. Si de plus $K/k$ est d'exposant fini, on v\'{e}rifie aussit\^{o}t que $B$ est une $r$-base de $K/k$ si et seulement si $B$ est un g\'{e}n\'{e}rateur minimal de $K/k$. Sous ces conditions, on d\'{e}signe par $di(K/k)=|G|$, o\`{u} $G$ est un g\'{e}n\'{e}rateur minimal de $K/k$, le degr\'{e} d'irrationalit\'{e} de $K/k$, et par $di(k)=|B|$, o\`{u} $B$ est une $r$-base de $k/k^p$, le degr\'{e} d'imperfection de $k$. Ces deux invariants permettent de mesurer respectivement  la taille de $K/k$ et la longueur de $k$.
Notamment, la taille d'une extension croit en fonction de l'inclusion. Plus pr\'{e}cis\'{e}ment, pour toute chaine d'extensions d'exposant born\'{e} $k\subseteq L\subseteq K$, on a $di(L/k)\leq di(K/k)$ (cf. \cite{Che1}). En particulier, cette propri\'et\'e permet d'\'etendre la meseure de la taille  \`a une extension purement ins\'{e}parable $K/k$ quelconque par prolongement vertical du degr\'{e} d'irrationalit\'{e} des sous-extensions interm\'{e}diaires d'exposant fini de $K/k$. Ainsi,  on pose $di(K/k)=\di\sup_{n\in\N}(di(k^{p^{-n}}\cap K /k))$, ici le $\sup$ est employ\'e dans le sens (\cite{N.B}, III, p. 25, proposition \textcolor{blue}{2}). Dans ce contexte, on montre dans (cf. \cite{Che1}) que la mesure de la taille d'une extension est compatible avec l'inclusion et la lin\'{e}arit\'{e} disjointe. En d'autres termes,  on a :
\sk

\begin{itemize}{\it
\item[$\bullet$] Pour toute chaine d'extensions purement ins\'{e}parables $k\subseteq L\subseteq L'\subseteq K$, on a $di(L'/L)\leq di(K/k)\leq di(k)$
\item[$\bullet$] Pour tous corps interm\'{e}diaires $K_1$ et $K_2$ de $K/k$, $k$-lin\'{e}airement disjoints, on a $di(K_1(K_2)/k)=di(K_1/k)+di(K_2/k)$ et $di(K_1(K_2)/K_2)=di(K_1/k)$.}
\end{itemize}
\sk

\noindent En outre, $di(K/k)=\sup(di(L/k))$, o\`{u} $k\subseteq L\subseteq K$. Autrement dit, la mesure de la taille de $K/k$ est vue comme limite inductive du degr\'e d'irrationalit\'e de ces sous-extensions interm\'{e}diaires.

Dans cette note, nous continuons \`{a} s'int\'{e}resser \`{a} la relation qui relit la r\'{e}partition restreinte de la finitude au niveau des corps interm\'{e}diaires d'une extension purement ins\'{e}parable $K/k$ \`{a} la finitude absolue ou globale de $K/k$.

Dans \cite{Dev2}, J. K. Deveney construit un exemple d'extension modulaire $K/k$ dite extension $w_0$-g\'en\'er\'ee  tel que toute sous-extension propre de $K/k$ est finie, et telle que pour tout $n\in \N$, on a $[k^{p^{-n}}\cap K: k]= p^{2n}$.  Il est facile de v\'{e}rifier que
$di(k^{p^{-n}}\cap K/k))=2$, et donc  $K/k$ est relativement parfaite d'exposant non born\'{e} dont la mesure de la taille vaut 2. En particulier, $K/k$ ne conserve pas la distribution de la finitude au niveau local des corps interm\'{e}diaires de $K/k$. Ainsi,  on peut esp\'{e}rer \'{e}tendre la $w_0$-g\'en\'eratrice \`{a} une extension de taille quelconque.
Dans  cette perspective, dans \cite{Che-Fli2}  on construit pour tout entier $j$ une extension purement ins\'{e}parable $K/k$ d'exposant non born\'{e}  v\'{e}rifiant :
\sk

\begin{itemize}{\it
\item[\rm(i)] Toute sous-extension propre de $K/k$ est finie;
\item[\rm(ii)] Pour tout $n\in \N$, $[k^{p^{-n}}\cap K: k]= p^{jn}$.}
\end{itemize}
\sk

\noindent Am\'{e}liorant ainsi le contre-exemple de J. K. Deveney, une telle extension est relativement parfaite d'exposant non born\'{e}, et pour tout $n\in \N$, $di(k^{p^{-n}}\cap K/ k)=j$. Il est \'{e}galement clair que  $K/k$ ne  conserve pas  la  finitude restreinte.
Il s'agit donc d'une forme d'irr\'{e}ductibilit\'{e} dans le sens o\`{u} $K/k$ ne peut se d\'{e}composer  sous la forme $k \longrightarrow K_1 \longrightarrow  K$ avec $K_1/k$  et $K/K_1$ ont chacune un exposant non born\'{e}. D'autre part, toute extension de taille finie est compos\'{e}e d'extensions $w_0$-g\'en\'er\'ees.
Toutefois, lorsque la taille de $K/k$ n'est pas born\'{e}, on montre que toute extension modulaire d'exposant non born\'{e} admet une sous-extension d'exposant non born\'{e}. En particulier, on montre pour qu'une extension $w_0$-g\'en\'er\'ee soit de taille finie il faut et il suffit que la plus petite sous-extenion $m$ de $K/k$ telle que $K/m$ est modulaire soit non triviale ($m\not= K$), et par suite si l'on tient compte de ce r\'esultat, il est fort probable que la $w_0$-g\'en\'eratrice soit li\'ee aux extensions de taille finie.  Ceci, nous am\`{e}ne \`{a} \'{e}tudier de pr\`{e}s la $w_0$-g\'en\'eratrice au sens restreint.  Conform\'{e}ment \`{a} cette approche,  et dans le but d'\'{e}tendre la $w_0$-g\'en\'eratrice aux extensions purement ins\'{e}parables de taille non born\'{e}e,  on propose d'autres g\'{e}n\'{e}ralisations. Une extension $K/k$ est dite $j$-$w_0$-g\'en\'er\'ee si $K/k$ n'admet aucun corps interm\'{e}diaire $L$ d'exposant non born\'{e} sur $k$ et de degr\'{e} d'irrationalit\'{e} inf\'{e}rieur ou \'{e}gal $j$. Il s'agit d'une forme d'irr\'edictubilit\'e locale conditionn\'ee par la mesure de la taille.   En particulier,  si pour tout entier $j$, $K/k$ est $j$-$w_0$-g\'en\'er\'ee,   $K/k$ sera appel\'{e}e $+\infty$-$w_0$-g\'en\'er\'ee. On v\'{e}rifie aussit\^{o}t que toute extension  $w_0$-g\'en\'er\'ee est $+\infty$-$w_0$-g\'en\'er\'ee et inversement toute extension $+\infty$-$w_0$-g\'en\'er\'ee de taille finie est $w_0$-g\'en\'er\'ee. Il s'agit d'une r\'{e}partition  absolue de la $w_0$-g\'en\'eratrice au niveau des extensions de taille finie.  Par ailleurs, pour des raisons de la non-contradiction, on construit un exemple d'extension $+\infty$-$w_0$-g\'en\'er\'ee de taille infinie.
\sk

Enfin, il est \`a noter qu'au cours de cette note, on reprend,  les notations et les r\'{e}sultats \'el\'ementaires de \cite{Che1},
puisqu'ils sont utilis\'{e}s avec toute leur force ici.
%%%%%%%%%%%%%%%%%%%%%%%%%%%%%%%%%%%%%%%
%%%%%%%%%%%%%%%%%%%%%%%%%%%%%%%%%%%%%%%%%%%%%%%%%%%%%%%%%%%%%%%%%%%%%%%%%%
\section{G\'{e}n\'{e}ralit\'{e}}

D'abord, nous commencerons par donner une  liste pr\'eliminaire  des notations le plus souvent utilis\'ees  tout le long de ce travail :
%on se servira des notations suivantes :
\sk

\begin{itemize}{\it
\item $k$ d\'esigne toujours un corps commutatif de caract\'eristique $p>0$, et $\Om$ une cl\^oture alg\'ebrique de $k$.
\item $k^{p^{-\infty}}$ indique la cl\^oture purement ins\'eparable de $\Om/k$.
\item Pour tout $a\in \Om$, pour tout $n\in \N^*$, on symbolise la racine du polyn\^ome $X^{p^n}-a$ dans $\Om$  par $a^{p^{-n}}$. En outre, on pose $k(a^{p^{-\infty}})=k(a^{p^{-1}},\dots, a^{p^{-n}},$ $\dots )=\di\bigcup_{n\in \N^*} k(a^{p^{-n}})$ et
$k^{p^{-n}}=\{a\in \Om\,|\ , a^{p^{n}}\in k\}$.
\item Pour toute famille $B=(a_i)_{i\in I}$  d'\'el\'ements de $\Om$, on note $k(B^{p^{-\infty}})= k({({a_i}^{p^{-\infty}})}_{i\in I})$.
\item Enfin, |.| sera employ\'e au lieu du terme cardinal.}
\end{itemize}
\sk

 Il est \`a signaler aussi que toutes les extensions qui interviennent dans ce papier sont des sous-extensions purement ins\'eparables de $\Om$, et il est commode de noter $[k,K]$ l'ensemble des corps interm\'{e}diaires d'une extension $K/k$.
\subsection{$r$-base, $r$-g\'en\'erateur}
\begin{df}Soit $K/k$ une extension. Une partie $G$ de $K$
est dite $r$-g\'en\'era\-t\-eur de $K/k$, si $K=k(G)$ ; et si de plus pour
tout $x\in G$, $x\not \in k(G\backslash x)$, $G$ sera appel\'{e}e
$r$-g\'en\'erateur minimal de $K/k$.
\end{df}
\begin{df}Etant donn\'ees  une extension $K/k$  de caract\'eristique $p>0$  et
une partie $B$ de $K$. On dit que $B$ est une  $r$-base de $K/k$, si $B$
est un $r$-g\'en\'erateur minimal de $K/k(K^p)$. Dans le m\^eme ordre d'id\'ees, on dit que $B$ est
$r$-libre sur $k$, si $B$ est une $r$-base de $k(B)/k$~; dans le cas
contraire $B$ est dite $r$-li\'ee sur $k$.
\end{df}

Voici quelques cas particuliers :
\sk

\begin{itemize}{\it
\item Toute $r$-base de $k/k^p$ s'appelle $p$-base de $k$.
\item Egalement, toute partie d'\'el\'ements de $k$, $r$-libre sur $k^p$ sera appel\'ee $p$-ind\'epend\-a\-n\-te (ou $p$-libre) sur $k^p$.}
\end{itemize}
\sk

Ici $B$ d\'esigne une partie d'un corps commutatif $k$ de caract\'eristique $p>0$. Comme cons\'equeces imm\'ediates on a :
\sk

\begin{itemize}{\it
\item[\rm{(1)}] $B$  est $p$-base de $k$ si et seulement si pour tout $n\in \Z$, $B^{p^n}$ l'est \'egalement  de $k^{p^n}$.
\item[\rm{(2)}] $B$ est $r$-libre sur $k^p$ si et seulement si pour tout $n\in \Z$, $B^{p^n}$ l'est auusi sur  $k^{p^{n+1}}$.
\item[\rm{(3)}] $B$  est $p$-base de $k$ si et seulement si $B$ est un $r$-g\'en\'erateur minimal de $k/k^p$.
\item[\rm{(4)}] $B$  est $p$-base de $k$ si et seulement si pour tout $n\in \N^*$, $k^{p^{-n}}=\otimes_k (\otimes_k $ $k($ $a^{p^{-n}}$ $))_{a\in B}$
et  pour tout $a\in B$, $a\not\in k^p$. En particulier, $B$  est $p$-base de $k$ si et seulement si $k^{p^{-\infty}}=\otimes_k (\otimes_k k(a^{p^{-\infty}}))_{a\in B}$ et pour tout $a\in B$, $a\not\in k^p$.}
\end{itemize}
\sk

Il est \`a noter que le produit tensoriel  est  utilis\'e conform\'ement  \`a la d\'efinition {5} (cf. \cite{N.B2}, III, p. 42). Il est vu comme limite inductive du produit tensoriel d'une famille finie de $k$-alg\`ebre.
Toutefois, la proposition ci-dessous permet de ramener  l'\'{e}tude des propri\'{e}t\'{e}s des syst\`{e}mes  $r$-libres des extensions  de haureur $\leq 1$,  ($K^p\subseteq k$) au cas fini. Plus pr\'{e}cis\'{e}ment, on a :
\begin{pro} {\label{pr1}} Soit $K/k$ une extension de caract\'{e}ristique $p>0$. Une partie $B$ de $K$ est $r$-libre sur $k(K^p)$ si et  seulement s'il en est de m\^{e}me pour toute sous-partie finie de $B$.
\end{pro}
\pre Imm\'{e}diat.  \cqfd
\begin{pro} {\label{pr2}} Soit $K/k$ une extension de caract\'{e}ristique $p>0$. Toute partie finie $B$ de $K$ satisfait  $[k(K^p)(B) :k(K^p)]\leq p^{|B|}$, et il y'a \'{e}galit\'{e} si et seulement si $B$ est $r$-libre sur $k(K^p)$.
\end{pro}
\pre  Notons $B=\{x_1,\dots,x_n\}$, comme pour tout $i\in \{1,\dots,n\}$, on a $x_{i}^p \in k(K^p)\subseteq k(K^p)(x_1,\dots,x_{i-1})$, alors $[k(K^p)(x_1,\dots,x_{i}) : k(K^p)(x_1,\dots,x_{i-1})]$ $\leq p$, et il y'a  \'{e}galit\'{e} si et seulement si $x_i\not\in k(K^p)(x_1,\dots,x_{i-1})$. Compte tenu de la transitivit\'e de la finitude, on a
$[k(K^p)(x_1,\dots,x_n) :k(K^p)]=\di\prod_{i=1}^{n}[k(K^p)(x_1,$ $\dots,$ $x_{i}) :k(K^p)(x_1,\dots,x_{i-1})]\leq p^n$, et il y'a \'{e}galit\'{e} si et seulement si $B$ est $r$-libre sur $k(K^p)$. \cqfd
\begin{cor} {\label{cor1}} Soit $K/k$ une extension de caract\'{e}ristique $p>0$. Une partie  $B$ de $K$ est $r$-libre sur $k(K^p)$ si et seulement si pour toute sous-partie finie $B'$ de $B$, on a $[k(K^p)(B') :k(K^p)]=p^{|B'|}$.
\end{cor}

Comme application, le r\'{e}sultat ci-dessous montre que la $r$-ind\'{e}pendance est transitive dans le cas des extensions de hauteurs $1$. Autrement dit :
\begin{pro} {\label{pr3}} Etant donn\'ee  une extension $K/k$ de caract\'eristique $p>0$. Deux parties
$B_1$ et $B_2$ de $K$ sont respectivement  $r$-libres sur $k(K^p)$ et $k(B_1)(K^p)$
si et seulement si  $B_1\cup B_2$ l'est
sur $k(K^p)$. En particulier, si $B_1$ est une $r$-base de
$k(K^p)/K^p$, et $B_2$ est une $r$-base de $K/k(B_1)(K^p)$, alors $B_1\cup
B_2$ est  $p$-base de $K$.
\end{pro}
\pre La condition suffisante r\'esulte aussit\^ot de la d\'efinition des $r$-bases. Par ailleurs, d'apr\`{e}s la proposition {\ref{pr1}}, on se ram\`{e}ne   au cas o\`{u} $B_1$ et $B_2$ sont finies.  En vertu de la proposition \textcolor{blue}{\ref{pr2}},  on a $[k(K^p)(B_1\cup B_2) :k(K^p)]= [k(K^p)(B_1\cup B_2) :k(K^p)(B_1)]. [k(K^p)(B_1) :k(K^p)]=p^{|B_2|}.p^{|B_1|}=p^{|B_2|+|B_1|}$ $=p^{|B_1\cup B_2|}$, et par suite $B_1\cup B_2$ est $r$-libre sur $k(K^p)$.  \cqfd \sk

Comme cons\'{e}quences imm\'ediates :
\begin{cor} {\label{adcor}} Soient $k\subseteq L\subseteq K$ des extensions purement ins\'eparables et $B_1$, $B_2$ deux parties respectivement de $K$ et $L$. Si $B_1$ est une $r$-base de $K/L$ et $B_2$ une $r$-base de $L(K^p)/k(K^p)$, alors $B_1\cup B_2$ est une $r$-base de $K/k$.
\end{cor}
\begin{cor} {\label{cor2}}  Soient $K/k$ une extension de caract\'eristique $p>0$, $x$ un \'el\'em\-ent de $K$, et $B$
une partie $r$-libre sur $k(K^p)$. Pour que $B\cup \{x\}$ soit  $r$-libre sur $k(K^p)$ il faut et il suffit que $x\not\in k(K^p)(B)$.
\end{cor}
\pre imm\'{e}diat. \cqfd
\begin{thm} {\label{thm1}} [th\'eor\`eme de la $r$-base incompl\`ete] Etant donn\'{e}es  une extension $K/k$ de caract\'eristique $p>0$, et une partie $B$ de $K$,  $r$-libre sur $k(K^p)$. Pour tout
$r$-g\'en\'erateur $G$ de $K/k(K^p)$, il existe un sous-ensemble $G_1$ de $G$
tel que $B\cup G_1$ est une $r$-base de $K/k$.
\end{thm}
\pre  Le cas o\`u $k(K^p)(B)=K$ est trivialement \'{e}vident.  Si $k(K^p)(B)\not= K$, il existe $x\in G$ tel que
$x\not\in k(K^p)(B)$. En effet, si pour tout $x\in G$, $x\in
k(K^p)(B)$,  comme $G$ est un $r$-g\'en\'erateur de $K/k(K^p)$, on aura $k(K^p)(G)=K\subseteq
k(K^p)(B)$, absurde. D'apr\`es le lemme pr\'ec\'edent, $B\cup \{x\}$
est une partie $r$-libre sur $k(K^p)$. Posons ensuite $H=\{L\subset G$ tel
que $B\cup L$ est $r$-libre sur $k(K^p)\}$. IL est clair que $H$ est inductif, et donc d'apr\`es
le lemme de Zorn, $H$ admet un \'el\'ement maximal que l'on note $M$. Soit  $B_1=M\cup B$, n\'ecessairement $K=k(K^p)(B_1)$, si
$K\not= k(K^p)(B_1)$, il existe \'{e}galement un \'el\'ement $y$ de $G$ tel que $y\not \in
k(K^p)(B_1)$, et donc $B_1\cup \{y\}$ serait
 $r$-libre sur $k(K^p)$~;  c'est une contradiction avec
le fait que $M$ est maximal.\cqfd \sk

Voici quelques cons\'equences imm\'ediates :
\sk

\begin{itemize}{\it
\item[\rm{(1)}]  De tout $r$-g\'en\'erateur de $K/k(K^p)$ on
peut en extraire  une $r$-base de $K/k$.
\item[\rm{(2)}]  Toute partie $r$-libre sur $k(K^p)$ peut
\^{e}tre compl\'et\'{e}e en une $r$-base de $K/k$. En particulier, toute partie $p$-ind\'ependante sur $k^p$ peut \^{e}tre \'{e}tendue en une $p$-base de $k$.
\item[\rm{(3)}]  Toute extension $K/k$ admet une $r$-base. En outre, tout corps commutatif de caract\'eristique $p>0$ admet une $p$-base.}
\end{itemize}
\sk

Par ailleurs, toutes les $r$-bases d'une m\^{e}me extension ont m\^{e}me cardinal comme le pr\'ecise le  r\'{e}sultat suivant.
\begin{thm} {\label{thm2}} Soit $K/k$ une extension de caract\'eristique $p>0$. Si $B_1$
et $B_2$  sont deux $r$-bases de $K/k$, alors $|B_1|=|B_2|$.
\end{thm}

Pour la preuve de ce th\'eor\`eme on se s\'ervira des r\'esultats
suivants.
\begin{lem} {\label{lem1}} {\bf [Lemme d'\'echange]} Sous les conditions du th\'eor\`eme
pr\'ec\'edent, pour tout $x\in B_2$, il existe $x_1\in B_1$ tel
que $(B_1\backslash \{x_1\})\cup \{x\}$ est une $r$-base de $K/k$.
\end{lem}
\pre  Choisissons un \'{e}l\'{e}ment arbitraire $x$ de $B_2$, comme $B_2$ est une $r$-base de $K/k$, il en r\'esulte que
 $\{x\}$ est $r$-libre sur $k(K^p)$. Compte tenu du
th\'eor\`eme \textcolor{blue}{\ref{thm1}}, il existe $B'_1\subset B_1$ tel que $B'_1\cup
\{x\}$ est une $r$-base de $K/k$.  D'o\`{u}, $p=[k(K^p)(B'_1)(\{x\}) :k(K^p)(B'_1)] =[K :k(K^p)(B'_1)]=[k(K^p)(B'_1)( B_1\backslash B'_1) :k(K^p)(B'_1)]$, et comme
 $B_1\backslash B'_1$ est $r$-libre sur  $k(K^p)(B'_1)$,   on en d\'{e}duit que $|B_1\backslash B'_1|=1$, c'est-\`{a}-dire $ B_1\backslash B'_1$ est r\'eduit \`a un singleton. \cqfd
\begin{pro} {\label{pr4}}Soit $K/k$ une extension de caract\'eristique $p>0$. Si $K/k$ admet au moins  une $r$-base finie, alors  toutes les $r$-bases de $K/k$ sont finies et ont m\^eme cardinal.
\end{pro}
\pre Imm\'ediat, il suffit d'appliquer la proposition \textcolor{blue}{\ref{pr2}}. \cqfd \sk

\noindent{} {\bf Preuve du th\'eor\`eme \textcolor{blue}{\ref{thm2}}.} D'apr\`es la prroposition \textcolor{blue}{\ref{pr1}}, on se ram\`ene  au cas  o\`u  $|B_1|$ et $|B_2|$ sont infinies.
Comme $B_1$ est une $r$-base de $K/k$, pour tout $x\in B_2$,  il existe une partie finie  $D(x)$ de $B_1$ telle
que $x\in k(K^p)(D(x))$, et par suite $K=k(K^p)(B_2)\subseteq k(K^p)(\di\bigcup_{x\in B_2}(D(x)))$.  Il en r\'esulte que $\di\bigcup_{x\in
B_2}(D(x))=B_1$, et  en vertu du (\cite{N.B}, III, p. 49, cor 3), on obtient $|B_1|\leq
|B_2|.|\N|=|B_2|$. De la m\^eme fa\c{c}on on montre que
$|B_2|\leq |B_1|$~;  d'o\`u
$|B_1|=|B_2|$. \cqfd \sk

Comme cons\'equence, on a :
\begin{cor} {\label{cor3}}
Pour toute  partie $B_1$ de $K$, $r$-libre sur $k(K^p)$, et tout
$r$-g\'en\'erateur $G$ de $K/k$, on a $|B_1|\leq |G|$.
\end{cor}
\pre Imm\'ediat, puisque tout $r$-g\'en\'erateur peut se r\'{e}duire (respectivement toute famille $r$-libre peut se compl\'eter) en une $r$-base. \cqfd \sk

Dans le cas o\`u $K/k(K^p)$ est finie, compte tenu du th\'eor\`eme de la $r$-base incompl\`ete, un $r$-g\'en\'erateur $G$ de $K/k(K^p)$ est une $r$-base de $K/k$ si et seulement si $|G|=Log_p([K:k(K^p)])$. En particulier, si $B$ est une $r$-base de $K/k$ et $G$ un
$r$-g\'en\'erateur de $K/k(K^p)$ tels que $|B|=|G|<+\infty$, alors $G$ est une $r$-base de $K/k$.
\sk

Soit $K /k$ une extension purement ins\'{e}parable de caract\'{e}ristique $p>0$. On rappelle que $K$ est dit d'exposant fini sur $k$,
s'il existe $e\in \N$ tel que $K^{p^e}\subseteq k$, et le plus petit entier qui satisfait cette relation sera appel\'e exposant (ou hauteur) de $K/k$. Certes,
la proposition  suivante  permet de ramener l'\'{e}tude des propri\'{e}t\'{e}s  des $r$-g\'en\'erateurs minimals des extensions
d'exposant fini au cas des extensions  de hauteur $1$,  lesquelles sont  plus riches.
\begin{pro} {\label{pr5}} Soit $K/k$ une extension purement ins\'eparable d'exposant fini. Pour qu'une  partie de $K$ soit une $r$-base de $K/k$ il faut et il suffit que elle soit $r$-g\'en\'erateur minimal de $K/k$.
\end{pro}
\pre
Soit $G$ une $r$-base de $K/k$, donc $K=k(K^p)(G)=\dots= k(K^{p^e})(G)$ $=k(G)$, et s'il existe $x\in G$ tel que $x\in
k(G\backslash \{x\})$, on aura $x\in k(K^p)(G\backslash \{x\})$, c'est une contradiction avec le fait que
$G$ est une $r$-base de $K/k$. Inversement, pour tout $r$-g\'en\'erateur minimal $G$  de $K/k$, on a
$K=k(G)=k(K^p)(G)$, et s'il existe $x\in G$ tel que $x\in K=
k(K^p)(G\backslash x)=\dots =k(K^{p^e}) (G\backslash
\{x\}=k(G\backslash \{x\})$, on aura une contradiction avec le fait que $G$ est
un $r$-g\'en\'erateur minimal de $K/k$. \cqfd
\begin{thm} {\label{thm3}} Soit $L/k$ une sous extension d'une extension purement ins\'e\-p\-a\-r\-a\-ble d'exposant fini $K/k$. Pour toutes $r$-bases $B_L$ et $B_K$ respectivement de
$L/k$ et $K/k$, on a $|B_L|\leq |B_K|$.
\end{thm}
\pre On distingue deux cas :
\sk

1-ier cas. Si $K/k$ est d'exposant 1, c'est-\`{a}-dire $K^p\subseteq k$, donc $L^p\subseteq k$. D'apr\`{e}s le th\'eor\`eme \textcolor{blue}{\ref{thm1}}, il existe $B_1\subseteq B_K$
tel que $B_L\cup B_1$ est une $r$-base de $K/k$, et par suite $|B_L|\leq |B_L\cup B_1|=|B_K|$.
\sk

2-i\`{e}me cas. Etant donn\'e  un entier naturel $e$ distinct de $0$ et $1$. Raisonnons par r\'ecurrence en supposant  que le th\'{e}or\`{e}me est v\'{e}rifi\'e  pour toute extension d'exposant $<e$, et soit $K/k$ une extension purement ins\'eprable d'exposant $e$.  Il est clair que  $k(K^p)\subseteq L(K^p)\subseteq K$, et donc il existe  $B_1\subseteq B_L$ et $B_2\subseteq B_K$ telles que $B_1$ et $B_2$ sont deux $r$-bases respectivement de $L(K^p)/k(K^p)$ et $K/L(K^p)$. D'apr\`es la transitivit\'e de la $r$-ind\'ependance, $B_1\cup B_2$ est une $r$-base de $K/k(K^p)$. Posons ensuite $k_1=k(B_1)$ et  $B'_L=B_L\setminus B_1$ ; on v\'erifie aussit\^ot que  $L\subseteq k_1(K^p)=k_1({B_2}^p)$, et $k_1(K^p)/k_1$ est d'exposant $<e$. Par application de la propri\'et\'e de r\'ecurrence et du corollaire \textcolor{blue}{\ref{cor3}},
on obtient $|B'_L|\leq  |{B_2}^p|=|B_2|$. Comme $B_1\cap B'_L=\emptyset$ et $B_1\cap B_2=\emptyset$, alors $|B_1\cup B'_L|\leq  |B_1\cup B_2|$, et par suite $|B_L|\leq |B_K|$. \cqfd
\section{Degr\'{e} d'irrationalit\'{e}}

 Soit $K/k$ une extension purement ins\'{e}parable. D\'esormais, et sauf mention expresse du contraire, pour tout $n\in \N^{*}$, on note $k_n=k^{p^{-n}}\cap K$, on obtient  ainsi  $k\subseteq k_1\subseteq \dots \subseteq k_n\subseteq \dots \subseteq K$, et $k_n/k$ est d'exposant fini. Soit $B_n$ une $r$-base de $k_n/k$, d'apr\`{e}s le th\'{e}or\`{e}me \textcolor{blue}{\ref{thm3}}, $|B_n|\leq |B_{n+1}|$.
 Ensuite, on pose  $di(K/k)= \di \sup_{n\in \N^*}(|B_n|)$, on rappelle que le $\sup$ est utilis\'e ici au sens du (\cite{N.B}, III, p. 25, proposition {2}).
\begin{df}
 L'invariant $di(K/k)$ d\'{e}fini ci-dessus s'appelle le degr\'{e} d'irrat\-io\-n\-al\-i\-t\'{e} de $K/k$.
\end{df}

En particulier, et pour des raisons de diff\'erenciation,  le degr\'{e} d'irrationalit\'{e} de $k/k^p$ sera appel\'e degr\'e d'imperfection de $k$ et sera not\'e $di(k)$.
\begin{rem}
$di(K/k)$ permet de mesurer la taille de l'extension $K /k$, et $di(k)$ la longueur  de $k$.
\end{rem}

Toutefois, on v\'erifie aussit\^ot que :
\sk

\begin{itemize}{\it
\item $di(K/K)=0$.
\item Pour tout $n\in \Z$, $di(k)=di(k^{p^n})=di(k^{p^{-\infty}}/k)$.
\item Compte tenu du corollaire \textcolor{blue}{\ref{adcor}}, pour toute sous-extension $L/k$ de $K/k$, on a $di(K/k(K^p))=di(K/L(K^p))+di(L(K^p)/k(K^p))$. Plus g\'en\'eralement, si $K/k$ est d'exposant $1$, on a $di(K/k)=di(K/L)+di(L/k)$.
\item En vertu de la proposition \textcolor{blue}{\ref{pr2}}, pour toute extension purement ins\'eparable d'exposant fini $K/k$, on a $di(K/k)=di(K/k(K^p))$.}
\end{itemize}
\begin{thm} {\label{thm4}} Soient $k\subseteq L\subseteq K$ des extensions purement ins\'{e}parables, on a
$di(L/k)\leq di(K/k)$. En outre, $di(K/k)=\di\sup(di(L/k))_{L\in [k, K]}$.
\end{thm}
\pre D'apr\`{e}s le th\'eor\`eme \textcolor{blue}{\ref{thm3}}, il suffit de remarquer  que pour tout $n\geq 1$,  on a $di(k^{p^{-n}}\cap L/k)\leq di(k_n/k)$, et donc $\di\sup(di(k^{p^{-n}}\cap L/k))_{n\geq 1}\leq \di\sup(di(k_n/k))_{n\geq 1}$ ; ou encore $di(L/k)\leq di(K /k)$. \cqfd \sk

Une cons\'{e}quence type est le r\'{e}sultat suivant :
\begin{thm} {\label{thm5}}
Pour toute extension purement ins\'{e}parable $K /k$, on a $di(K /k)$ $\leq di(k)$.
\end{thm}
\pre Il  suffit de remarquer qu'une partie $B$ de $k$ est une $p$-base de $k$ si et seulement si  $B^{p^{-n}}$ est  une $r$-base de $k^{p^{-n}}/k$ pour tout $n\geq 1$. Comme $k^{p^{-\infty}}=\di\bigcup_{n\geq 1} k^{p^{-n}}$, on a  pour tout $n\geq 1$, $k^ {p^{-n}}\cap K\subseteq k^{p^{-\infty}},$  et par suite $di(K/k)\leq di(k^{p^{-\infty}}/k)=di(k)$. \cqfd
\begin{pro} {\label{pr6}} Soit $(K_n/k)_{n\in \N}$ une famille croissante de sous-extensions purement ins\'{e}parables d'une extension $\Om/k$. On a : $$di(\di\bigcup_{n\in\N}(K_n)/k)= \di\sup_{n\in\N}(di(K_n/k)).$$
\end{pro}
\pre  Notons $K=\di\bigcup_{n\in\N}K_n$, et soit $j$ un entier naturel non nul. Il est imm\'ediat que $k_j=k^{p^{-j}}\cap K=\di\bigcup_{n\in \N} (k^{p^{-j}}\cap K_n)$. Dans la suite on distingue deux cas :
\sk

 1-ier cas : si $di(k_j/k)$ est fini, ou encore $k_j/k$ est finie. Comme pour tout $n\in N$, on a $k^{p^{-j}}\cap K_n\subseteq k^{p^{-j}}\cap K_{n+1}\subseteq k^{p^{-j}}\cap K$, alors la suite d'entiers $([k^{p^{-j}}\cap K_n:k])_{n\in \N}$ est croissante et born\'ee, donc stationnaire \`a partir d'un rang $n_0$ ; et par cons\'equent pour tout $n\geq n_0$, $k^{p^{-j}}\cap K_n=k^{p^{-j}}\cap K_{n+1}$. En outre, $di(k^{p^{-j}}\cap K/k)=di(k^{p^{-j}}\cap K_{n_0}/k)=\di\sup_{n\in N} (di(k^{p^{-j}}\cap K_n/k))$.
\sk

 2-i\`eme cas : si $di(k^{p^{-j}}\cap K/k)$ est infini, ou encore $\di\sup_{n\in N} (di(k^{p^{-j}}\cap K_n/k))$ n'est pas fini. Comme $k^{p^{-j}}\cap K=\di\bigcup_{n\in \N} (k^{p^{-j}}\cap K_n)$, donc si $B^{j}_{n}$ est une $r$-base de $k^{p^{-j}}\cap K_n/k$, alors $\di\bigcup_{n\in\N} B^{j}_{n}$ est un $r$-g\'{e}n\'{e}rateur de $k^{p^{-j}}\cap K/k$. En vertu du corollaire \textcolor{blue}{\ref{cor3}},  $di(k^{p^{-j}}\cap K/k)\leq |\di\bigcup_{n\in\N} B^{j}_{n}|$, et d'apr\`es (\cite{N.B}, III,  p.49, corollaire {3}),  $|\di\bigcup_{n\in\N} B^{j}_{n}|\leq \di\sup_{n\in\N}(|B^{j}_{n}|)=\di\sup_{n\in\N}(di(k^{p^{-j}}\cap K_n/k))$.
 \sk

 Compte tenu de ces deux cas, on en d\'eduit que  $di(K/k)\leq \di\sup_{n\in\N}(di(K_n/k))$. Mais comme $K_n\subseteq K$ pour tout $n\geq 1$, d'apr\`{e}s le th\'{e}or\'{e}me \textcolor{blue}{\ref{thm4}} on obtient $\di\sup_{n\in\N}(di(K_n/k))\leq di(K/k)$, et par suite $di(K/k)=\di\sup_{n\in\N}(di(K_n/k))$.\cqfd \sk

Le r\'{e}sultat  suivant qui est une cons\'{e}quence bien connue de la lin\'{e}arit\'{e} disjointe   intervient souvent dans le reste de ce papier.
\begin{pro} {\label{pr7}} Soient $K_1/k$ et $K_2/k$ deux sous-extensions d'une m\^eme extension $K/k$, $k$-lin\'{e}airement disjointes.  Pour touts corps interm\'{e}diaires  $L_1$ et $L_2$  respectivement de $K_1$ et $K_2$, on a $L_2(K_1)$ et $L_1(K_1)$ sont $k(L_1, L_2)$-lin\'{e}airement-disjointes. En particulier, $L_2(K_1)\cap L_1(K_2)=k(L_1,L_2)$.
\end{pro}

Une famille $(F_i/k)_{i\in J}$ d'extensions est dites $k$-lin\'{e}airement disjointes,  si pour toute partie  $G$ d'\'{e}l\'{e}ments finis de $J$, $(F_n/k)_{n\in G}$ sont $k$-lin\'{e}airement disjointes (cf. \cite{F-J}, p. 36). Il est trivialement \'{e}vident que $k((F_i)_{i\in J})=\di \prod_{i\in J} F_i\simeq \otimes_k (\otimes_k F_i)_{i\in J}$ si et seulement si $(F_i/k)_{i\in J}$ sont $k$-lin\'{e}airement disjointes. De plus, les propri\'{e}t\'{e}s de la lin\'{e}arit\'{e} disjointe du cas fini se prolonge naturellement \`{a} une famille quelconques d'extensions  $k$-lin\'{e}airement disjointes. En particulier,  pour tout $i\in J$, soit $L_i$ un sous-corps de $F_i$, si $(F_i/k)_{i\in J}$ sont $k$-lin\'{e}airement disjointes,  compte tenu de la transitivit\'{e} de la lin\'{e}arit\'{e} disjointe,  $(L_i/k)_{i\in J}$ (resp. $((\di \prod_{n\in J}L_n)F_i/k)_{i\in J}$) sont $k$-lin\'{e}airement (resp. $\di \prod_{n\in J}L_n$-lin\'{e}airement) disjointes.
\sk

Consid\'erons maintenant  deux sous-extensions $K_1/k$ et $K_2/k$ d'exposant fini d'une m\^eme extension purement ins\'eparable  $K/k$. On v\'erifie aussit\^ot que si $B_1$ et $B_2$ sont deux $r$-bases respectivement de $K_1/k$ et $K_2/k$, alors $B_1$ et $B_1\cup B_2$ sont  deux $r$-g\'en\'erateurs respectivement de $K_1(K_2)/K_2$ et $K_1(K_2)/k$. En outre, $di(K_1(K_2)/K_2)\leq di(K_1/k)$ et $di(K_1(K_2)/k)\leq di(K_1/k)+di(K_2/k)$. D'une fa\c{c}on  plus pr\'ecise, on a :
\begin{pro} {\label{pr8}}Sous les conditions ci-dessus, et si de plus $K_1/k$ et $K_2/k$ sont
$k$-lin\'{e}airement disjointes, on a :
\sk

\begin{itemize}{\it
\item[\rm{(i)}] $B_1\cup B_2$ est une $r$-base de $K_1(K_2)/k$.
\item[\rm{(ii)}] $B_1$ est une $r$-base de $K_1(K_2)/K_2$.}
\end{itemize}
\end{pro}
\pre Ici, on se contente de pr\'{e}senter uniquement la preuve du premier item, puisque  les deux assertions  utilisent les m\^{e}mes techniques de raisonnement.
Il est clair que $K_1(K_2)=k(B_1\cup B_2)$, il suffit donc de montrer que $B_1\cup B_2$ est  minimal.  Pour cela, on suppose par exemple l'existence d'un  \'el\'ement $x$ dans $B_1$ tel que $x\in k((B_1\setminus\{x\})\cup B_2)=K$.  Comme $K_1/k$ et $K_2/k$ sont $k$-lin\'{e}airement disjointes,  par transitivit\'{e}, on a $k(B_1)=K_1$ et $K_2(B_1\setminus \{x\})=K$ sont $k(B_1\setminus \{x\})$-lin\'{e}airement disjoints, et
donc $K_1=K\cap K_1=k(B_1\setminus \{x\})$, c'est une contradiction avec le fait que $B_1$ est une $r$-base de $K_1/k$. \cqfd \sk

Comme cons\'{e}quence imm\'{e}diate, on a
\begin{cor} {\label{cor4}}  Soient $K_1$ et $K_2$ deux corps interm\'{e}diaires d'une m\^{e}me extension purement ins\'{e}parable $\Om/k$. Alors :
\sk

\begin{itemize}{\it
\item[\rm{(i)}] $di(K_1(K_2)/k)\leq di(K_1/k)+di(K_2 /k)$, et il y'a \'egalit\'e si $K_1$ et $K_2$ sont $k$-lin\'{e}airement disjoints.
\item[\rm{(ii)}] $di(K_1(K_2)/K_2)\leq di(K_1/k)$, et il y'a \'egalit\'e si $K_1$ et $K_2$ sont $k$-lin\'{e}air\-e\-ment disjoints.}
\end{itemize}
\end{cor}
 \pre Il suffit de remarquer que $K_1(K_2)=\di\bigcup_{j\in \N} (k^{p^{-j}}\cap K_1)(k^{p^{-j}}\cap K_2)=\di\bigcup_{j\in \N} K_2(k^{p^{-j}}\cap K_1)$, et si $K_1$ et $K_2$ sont $k$-lin\'{e}airement disjoints, d'apr\`{e}s la transitivit\'{e} de la lin\'{e}arit\'{e} disjoint, $k^{p^{-j}}\cap K_1$ et $k^{p^{-j}}\cap K_2$ sont aussi $k$-lin\'{e}airement disjoints pour tout $j\geq 1$.  On se ram\`{e}ne ainsi  au cas o\`{u} $K_1/k$ et $K_2/k$  sont d'exposant fini auquel cas le r\'{e}sultat d\'{e}coule imm\'{e}diatement de la proposition pr\'ec\'edente. \cqfd \sk

 Comme cons\'equence imm\'ediate, on a :
 \begin{cor} {\label{cor5}} Pour toute sous-extension $L/k$ d'une extension purement ins\'e\-p\-a\-rable $K/k$, on a $di(L(K^p)$ $ /k(K^p))\leq di(L/k(L^p))$, et il y'a \'egalit\'e si  $k(K^p)$ et $L$ sont $k(L^p)$-lin\'eairement disjointes.
 \end{cor}
\pre Due au corollaire \textcolor{blue}{\ref{cor4}}. \cqfd \sk

Le r\'esultat suivant am\'eliore naturelement les conditions du th\'eor\`eme \textcolor{blue}{\ref{thm4}}
\begin{thm} {\label{thm6}} Pour toute famille  d'extensions purement ins\'eparables $k\subseteq L\subseteq L'\subseteq K$, on a $di(L/L')\leq di(K/k)$.
\end{thm}
\pre Il est clair que $K=\di\bigcup_{j\in \N} L(k_j)$, et d'apr\`es la proposition \textcolor{blue}{\ref{pr6}}, et le th\'eor\`eme \textcolor{blue}{\ref{thm4}}, on a $di(L'/L)\leq di (K/L)=\di\sup_{j\in \N}(di(L(k_j)/k))\leq \di\sup_{j\in \N}(di(k_j/k))=di(K/k)$. \cqfd \sk

Comme cons\'equence imm\'ediate, on a :
\begin{cor} {\label{cor6}} Pour toute extension purement ins\'eparable $K/k$, on a $di(K)\leq di(k)$.
\end{cor}
\pre Il suffit de remarque que $K\subseteq k^{p^{-\infty}}$, et $di(K)=di(K/K^p)\leq di(k^{p^{-\infty}}/k^p)=di(k)$.\cqfd
 %%%%%%%%%%%%%%%%%%%%%%%%%%%%%%%%%%%%%%%%µµµµµµµµµ
 \subsection{Extensions relativement parfaites}

Au cours de cette section, on reprend, en les am\'eliorant,  quelques notions et r\'esultats de {\cite{Che-Fli4}}, puisqu'ils sont utilis\'es fr\'equemment ici.
\sk

Un corps $k$ de caract\'eristique $p$ est dit parfait si $k^{p}=k$ ; dans le m\^{e}me ordre d'id\'{e}es,
on dit que $K/k$ est relativement parfaite  si $k(K^{p})=K$. On v\'erifie ais\'ement que :
\sk

\begin{itemize}{\it
\item La relation "\^etre relativement parfaite" est transitive, c'est-\`a-dire si $K/L$ et $L/k$ sont relativement parfaites, alors $K/k$ l'est aussi.
\item Si $K/k$ est relativement parfaite, il en est de m\^eme de $L(K)/k(L)$.
\item La propri\'et\'e "\^etre relativement parfaite" est stable par un produit quelconque portant sur $k$. Autrement dit, pour toute famille $(K_i/k)_{i\in I}$ d'ext\-e\-nsions relativement parfaites, on a alors $\displaystyle \di \prod_{i}^{}K_{i}/k$ est aussi relativement parfaite.}
\end{itemize}
\sk

\noindent Par suite, il existe une plus grande sous-extension relativement parfaite de $K/k$ appel\'ee cl\^oture relativement parfaite de $K/k$, et se note $rp(K/k)$.
On a les relations
d'associativit\'{e}-transitivit\'{e} suivantes.
\begin{pro} Soit $L$ un corps interm\'{e}diaire de
$K/k$. Alors
$$
rp(rp(K/L)/k)=rp(K/k) \quad \mbox{ et } \quad
rp(K/rp(L/k))=rp(K/k).
$$
\end{pro}
\pre Cf. {\cite{Che-Fli4}}, p. {50},  proposition {5.2}. \cqfd
\begin{cor}Pour tout $L\in
[k,K]$, on a $K/L \hbox{ finie} \Longrightarrow rp(K/k) \subset  L.$
\end{cor}

En particulier, si $K/k$ est relativement parfaite, on a $K/L \hbox{
$finie$ } \Longrightarrow   L=K.$
Sch\'{e}matiquement on a un $trou$
\sk

\[\begin{array}{rcl}
k\longrightarrow &&K;\\
&\uparrow&\\
&\hbox{$trou$ }&
\end{array}\]

\noindent et ce $trou$ caract\'{e}rise le fait que $K/k$ est
relativement parfaite. En effet,  supposons que
$K/k$ v\'{e}rifie le $trou$ et soit $B$ une $r$-base de $K/k$.
Supposons $B\neq \emptyset$;
soit $x\in B$ et $L=k(K^{p})(B\setminus\left \{x\right \})$; on a $K/L$ est
finie, donc $K=L$ ce qui est absurde.

\begin{pro} {\label{arpaa1}} Soit $K/k$ une extension purement ins\'eparable telle que $[K:k(K^p)]$ est fini. Alors on a :
\sk

\begin{itemize}{\it
\item[{\rm (i)}] $K$ est relativement parfaite sur une extension finie de $k$.
\item[{\rm (ii)}] La suite d\'{e}croissante  $(k(K^{p^{n}}))_{n \in
\N}$ est stationnaire sur $k(K^{p^{n_{0}}})=rp(K$ $/$ $k)$.}
\end{itemize}
\end{pro}
\pre Cf. {\cite{Che-Fli4}}, p.  {51},
 lemme {2.1}. \cqfd \sk

Comme  cons\'equence de la proposition pr\'ec\'edente, on a :
\begin{pro} Soit $K/k$ une extension purement ins\'eparable telle que $[K:k(K^p)]$ est fini. Pour tout $L\in [k,K]$, on a $rp(K/L)=L(rp(K/k)).$
\end{pro}
\pre Cf. {\cite{Che-Fli4}}, p. {51},
 proposition {6.2}. \cqfd \sk

 En utilisant le lemme 1.16 qui se trouve dans (\cite{Mor-Vin}, p. 10), on peut affirmer que la condition de finitude de $[K:k(K^p]$ est n\'ec\'essaire, et par suite, le r\'esultat pr\'ec\'edent peut tomber en d\'efaut si $K/k(K^p)$ n'est pas finie.
Par ailleurs, on v\'erifie aussit\^ot que $k(K^p)=rp(K/k)(K^p)$, et donc pour qu'une partie $G$ de $K$ soit $r$-base de $K/k$ il faut et il suffit qu'elle en soit de m\^eme de $K/rp(K/k)$. De plus, comme $2$-i\`eme cons\'equence de la proposition \textcolor{blue}{\ref{arpaa1}}, le r\'esultat suivant exprime une condition n\'ecessaire et suffisant pour que $K/rp(K/k)$ soit finie. Plus pr\'ecis\'ement, on a :
\begin{pro} {\label{arp1}} Soit $K/k$ une extension purement ins\'eparable, alors $K/rp($ $K/$ $k)$ est finie si est seulement il en est de m\^eme de $K/k(K^p)$.
\end{pro}
\pre R\'esulte de la proposition \textcolor{blue}{\ref{arpaa1}}. \cqfd
%%%%%%%%%%%%%%%%%%%%%%%%%%%%%%%%%%%%%%%%%%%%%%%%%%%%%%%%%%%%%%%%%%%%%%%%%%%%%%%%%%%%%%%%%
%%%%%%%%%%%%%%%%%%%%%%%%%%%%%%%%%%%%%%%%%%%%%%%%%%%%%%%%%%%%%%%%%%%%%%%%%%%%%%%%%%%%%%%%%
\section{Extensions $q$-finies}
\begin{df}  Toute extension de degr\'e d'irrationalit\'{e}  fini s'appelle extension $q$-finie.
\end{df}

En d'autres sens, la $q$-finitude est synonyme de la finitude horizontale. Toutefois, la finitude se traduit par la finitude horizontale et verticale, il s'agit de la finitude au point de vue taille et hauteur. Autrement dit, $K/k$ est finie si et seulement si $K/k$ est $q$-finie d'exposant born\'e. Par ailleurs, on v\'{e}rifie que {\it
le degr\'e d'irrationalit\'e d'une extension $K/k$ vaut 1 si est seulement si l'ensemble de corps interm\'ediaires de $K/k$ est totalement ordonn\'e.}
Ensuite, on appelle extension $q$-simple toute extension qui satisfait l'affirmation pr\'ec\'edente.
\begin{rem} On rappelle que lorsque $di(k)$ est fini, et apr\`es avoir montr\'e dans \cite{Che-Fli2} que $K/k(K^p)$ est finie et $di(K)\leq di(k)$,  le degr\'e d'irrationalit\'e d'une extension purement ins\'eparable $K/k$ a \'et\'e d\'efini  par l'entier  $di(K/k)=di(k)-di(K)+di(K/k(K^p))$. En outre, toute extension est $q$-finie si $di(k)$ est fini.  Avec quelques  modifications l\'eg\`eres, on peut toujours prolonger cette d\'efinition au cas o\`u $di(k)$ est non born\'e.  Commen\c{c}ons par le choix d'une extension $K/k$ relativement parfaite et $q$-finie.
Etant donn\'ee une $p$-base $B$ de $k$, donc $k=k^p(B)$, et par suite $k(K^p)=K^p(B)$. Comme $K/k$ est relativement parfaite, alors $K=k(K^p)=K^p(B)$. D'apr\`es le th\'eor\`eme \textcolor{blue}{\ref{thm1}}, il existe $B_1\subseteq B$ telle que $B_1$ est une $p$-base de $K$. Ainsi,  on aura $k^{p^{-\infty}}=k(B^{p^{-\infty}})=k({B_1}^{p^{-\infty}})\otimes_k  k((B\setminus B_1)^{p^{-\infty}})\simeq K^{p^{-\infty}}\simeq  K \otimes_k k({B_1}^{p^{-\infty}} )$. En particulier, d'apr\`es le corollaire \textcolor{blue}{\ref{cor4}}, $di(K/k)=di(K \otimes_k k({B_1}^{p^{-\infty}} )/ k({B_1}^{p^{-\infty}} ))=di(k^{p^{-\infty}}/ k({B_1}^{p^{-\infty}} ))=di(k((B\setminus B_1)^{p^{-\infty}})/k)=| B\setminus B_1|$. Si on interpr\`ete (par abus de langage) $| B\setminus B_1|$ comme diff\'erence de degr\'e d'imperfection de $k$ et $K$ en \'ecrivant $| B\setminus B_1|=di(k)-di(K)$, on obtiendra  $di(K/k)=di(k)-di(K)$.    Dans le cas g\'en\'eral, supposons que $K/k$ est $q$-finie quelconque, donc $K/rp(K/k)$ est finie, d'o\`u $di(K)=di(rp(K/k))$; et par suite  $di(K/k)=di(rp(K/k)/k)+di(K/k(K^p))=di(k)-di(K)+di(K/k(K^p))$ (cf. proposition \textcolor{blue}{\ref{pr11}} ci-dessous).
\sk

Il est \`a signaler en tenant compte de cette consid\'eration que tous les r\'esult\-a\-ts des articles \cite{Che-Fli1}, \cite{Che-Fli2}, \cite{Che-Fli3}, \cite{Che-Fli5} se g\'en\'eralisent naturellement par translation  \`a une extension $q$-finie quelconque.
\end{rem}

Soient $L/k$ une sous-extension d'une extension $q$-finie $K/k$, pour tout $n\in \N$, on note toujours $k_n=k^{p^{-n}}\cap K$. On v\'{e}rifie aussit\^{o}t que :
\sk

\begin{itemize}{\it
\item[\rm{(i)}] La $q$-finitude est transitive, en particulier, pour tout $n\in \N$, $K/k(K^{p^n})$ et $k_n/k$ sont finies.
\item[\rm{(ii)}] Il existe $n_0\in \N$, pour tout $n\geq n_0$, $di(k_n/k)=di(K/k)$.}
\end{itemize}
\sk

Par ailleurs, voici quelques applications imm\'ediates des propositions \textcolor{blue}{\ref{arpaa1}} et \textcolor{blue}{\ref{arp1}}.
\begin{pro} {\label{pr9}} Soit $K/k$ une extension $q$-finie. La suite $(k(K^{p^n}))_{n\in\N}$ s'arr\-\^{e}\-te sur $rp(K/k)$ \`{a} partir d'un $n_0$. En particulier,  $K/rp(K/k)$ est finie.
\end{pro}

Comme cons\'{e}quence, on a :
\begin{cor} {\label{cor7}} La cl\^{o}ture relativement parfaite d'une extension $q$-finie $K/k$ n'est pas triviale. Plus pr\'{e}cis\'{e}ment,  $rp(K/k)/k$ est d'exposant non born\'{e} si  $K/k$ l'est.
\end{cor}
\pre Imm\'ediat. \cqfd
\begin{pro} {\label{pr10}} Pour toute  extension $q$-finie $K/k$, il existe $n\in \N$ tel que $K/k_n$ est relativement parfaite. En outre,  $k_n(rp(K/k))=K$.
\end{pro}
\pre Imm\'ediat. \cqfd
\begin{pro} {\label{pr11}} Le degr\'{e} d'irrationalit\'{e} d'une extension $q$-finie  $K/k$ v\'{e}rifie l'\'{e}galit\'{e} suivante : $di(K/k)$ $=di(rp(K/k)/k)+di(K/k(K^p))=di(K/rp(K/k))$ $+di(rp(K/k)/k)$.
\end{pro}
\pre Soient $G$ une $r$-base de $K/k$ et $K_r=rp(K/k)$, donc $k(G)/k$ admet un exposant  fini not\'e $m$ et, $K=K_r(G)$. En paticulier,  pour tout $n\geq m$,  $k(G)\subseteq k_n$. Compte tenu de la $r$-ind\'{e}pendance de $G$ sur $k(K^p)$ et vu que $k({k_n}^p)$ est un sous-ensemble de $k(K^p)$, on en d\'{e}duit que $G$ est $r$-libre sur $k({k_n}^p)$ pour tout $n\geq m$. Compl\'{e}tons $G$ en une $r$-base de $k_n/k$ par une partie $G_n$ de $k_n$. Dans ces conditions, pour $n$ suffisamment grand, on aura  $|G|+|G_n|=\di \sup_{j\geq m}(|G|+|G_j|)=di(K/k)=di(K_r(G)/k) \leq di(K_r/k)+di(k(G)/k)=di(K_r/k)+|G|$, et donc $|G_n|\leq di(K_r/k)$. Toutefois, comme $\di\bigcup_{n\geq m} k({k_n}^{p^m})=\di\bigcup_{n\geq m} k({G_n}^{p^m},G^{p^m})=\di\bigcup_{n\geq m} k({G_n}^{p^m})=k(K^{p^m})=K_r(K^{p^m})$, d'apr\`{e}s le th\'{e}or\`{e}me \textcolor{blue}{\ref{thm4}}, pour $n$ suffisamment grand, on aura \'egalement $ di(K_r/$ $k)$ $\leq di(K_r(K^{p^{m}})/k)= di(k({k_n}^{p^m})/k)\leq |{G_n}^{p^m}|=|G_n|$. D'o\`u, $|G_n|$ $=di(K_r/k)$ pour $n$ assez grand, et par suite $di(K/k)$ $=di(K_r/k)+di(K/k(K^p))$. \cqfd
\sk

Comme cons\'equence imm\'ediate, on a :
\begin{cor} {\label{cor8}}Pour qu'une extension $q$-finie $K/k$ soit finie il faut et il suffit que $di(K/k)=di(K/k(K^p))$.
\end{cor}
\begin{thm} {\label{thm7}} Pour toutes extensions $q$-finies $k\subseteq L\subseteq K$, on a $di(K/k)\leq di(K/L)+di(L/k)$, avec l'\'{e}galit\'{e} si et seulement si $L/k(L^p)$ et $k(K^p)/k(L^p)$ sont $k(L^p)$-lin\'{e}airement disjointes.
\end{thm}
\pre Comme $K=\di\bigcup_{n\in \N} L^{p^{-n}}\cap K$ et $K/k$ est $q$-finie, d'apr\`{e}s le th\'{e}or\`{e}me \textcolor{blue}{\ref{thm4}}, pour $n$ assez grand, on a $di(K/k)=di(L^{p^{-n}}\cap K/k)$ ; donc on est amen\'e au cas o\`{u} $K/L$ est finie, ou encore $rp(K/k)=rp(L/k)$.  Dans la suite,  on posera $L_r=K_r=rp(K/k)$. D'apr\`{e}s la proposition \textcolor{blue}{\ref{pr11}} ci-dessus,  on aura $di(K/k)=di(K_r/k)+di(K/k(K^p))= di(L_r/k)+di(K/L(K^p))+di(L(K^p)/k(K^p))= di(L_r/k)+di(K/L)+di(L(K^p)/k(K^p))$. Compte tenu du corollaire \textcolor{blue}{\ref{cor5}}, on aura $di(L($ $K^p$ $)$ $/k(K^p))\leq di(L/k(L^p))$, et donc $$di(K/k)\leq di(L_r/k)+di(K/L) +di(L/k(L^p))=di(L/k)+di(K/L),$$ toutefois il y'a  \'{e}galit\'{e} si et seulement si $di(L/k(L^p))= di(L(K^p)/k($ $K^p))$, ou encore $[L : k(L^p)]=[L(K^p) :k(K^p)]$, c'est-\`a-dire $L/k(L^p)$ et $k(K^p)/k($ $L^p)$ sont $k(L^p)$-lin\'{e}airement disjointes.\cqfd
\begin{rem} La condition de la lin\'{e}arit\'{e} disjointe qui figure dans la proposition ci-dessus  se traduit en terme de $r$-ind\'{e}pendance par toute $r$-base de $L/k$ se compl\`{e}te en une $r$-base de $K/k$.
\end{rem}

Comme application imm\'{e}diate, on a :
\begin{cor}{\label{ccor1}} Toute sous-extension relativement parfaite $L/k$ d'une extension $q$-finie $K/k$ v\'erifie $di(K/k)=di(K/L)+di(L/k)$.
\end{cor}

D'une fa\c{c}on assez g\'en\'erale, on a :
\begin{pro} {\label{pr12}} Pour toute suite de sous-extensions relativement parfaites  $k=K_0 \subseteq K_1\subseteq \dots\subseteq K_n$ d'une extension $q$-finie $K/k$,  on a $di(K/k)=\di \sum_{i=0}^{n-1} di(K_{n+1}/K_n) +di(K/K_n)$.
\end{pro}
\pre R\'{e}sulte imm\'{e}diatement du corollaire pr\'{e}c\'{e}dent. \cqfd \sk

Dans la suite on va \'{e}tudier de plus pr\`{e}s  les propri\'{e}t\'{e}s des  exposants d'une extension $q$-finie.
\subsection{Exposants d'une extension $q$-finie}

Dans cette section nous distinguons  deux cas~:
\sk

{\bf Cas o\`u ${K/k}$ est purement ins\'{e}parable finie.}
Soit $x\in K$, posons $o( x/k ) = \inf\{$ $m \in \mathbf{N}|\; x^{p^m}\in k \}$
et $o_1(K/k) = \inf\{m\in \mathbf{N}|\; K^{p^{m}}\subset k\}$. Une
$r$-base $B=\{a_{1},a_{2},\dots, a_{n}\}$ de $K/k$ est dite
canoniquement ordonn\'{e}e si pour $j=1,2,\dots,n$, on a
$o(a_{j}/k(a_{1},a_{2},$ $\dots$ $ ,a_{j-1}))=
o_{1}(K/k(a_{1},a_{2},\dots ,a_{j-1})).$ Ainsi, l'entier
$o(a_j/k(a_1,\dots, $ $a_{j-1}))$ d\'{e}fini ci-dessus v\'{e}rifie
$o(a_j/k(a_1,\dots, a_{j-1}))=\inf\{m$ $\in \mathbf{N}|
\; di(k(K^{p^m})/k)\leq j-1\}$ (cf. {\cite{Che-Fli2}}, p.
{138}, lemme {1.3}). On en
d\'eduit aussit\^ot  le r\'{e}sultat de ({\cite{Pic}}, p.
{90}, satz {14}) qui confirme
l'ind\'{e}p\-e\-n\-d\-a\-n\-ce des entiers $o(a_i/k(a_1,$ $\dots ,
a_{i-1}) )$, $(1\leq i\leq n)$, vis-\`a-vis au choix des $r$-bases
canoniquement ordonn\'{e}es $\{a_1,\dots , $ $a_n\}$  de $K/k$. Par
suite, on pose  $o_i(K/k)=o(a_i/k(a_1,$ $\dots , a_{i-1}) )$ si
$1\leq i\leq n$, et $o_i(K/k)=0$ si $i>n$, o\`u $\{a_1,\dots , a_n\}$
est une $r$-base canoniquement ordonn\'{e}e de $K/k$. L'invariant $o_i(K/k)$
ci-dessus s'appelle le $i$-\`{e}me exposant de $K/k$. Voici les principales relations dont on aura besoin, et qui
font intervenir les exposants.
\begin{pro} {\label{pr13}}
Soient $K$ et $L$ deux corps interm\'{e}diaires d'une
extension $\Omega /k$, avec $K/k$ purement ins\'{e}parable finie.
Alors pour tout entier $j$, on a $o_{j}(K($ $L)/k(L))\leq
o_{j}(K/k)$.
\end{pro}
\pre Cf. {\cite{Che-Fli}}, p. {373},
 proposition {5}. \cqfd
\begin{pro} {\label{pr14}}
Soit $K/k$ une extension purement ins\'{e}parable
finie. Pour toute sous-extension $L/L'$ de $K/k$, et  pour tout
$j\in\mathbf{N}$, on a  $o_j(L/L')\leq o_j(K/k)$.
\end{pro}
\pre cf. {\cite{Che-Fli}},
p. {374}, proposition {6}. \cqfd
\begin{pro} {\label{pr15}}
 Soient $\{\al_1 ,\dots,\al_n\}$
une $r$-base canoniquement ordonn\'{e}e de $K/k$, et $m_j$ le $j$-i\`{e}me exposant de $K/k$,  $1\leq j\leq n$. On a~:
\sk

\begin{itemize}{\it
\item[{\rm(1)}]  $k(K^{p^{m_j}})=k(\al_{1}^{p^{m_j}},\dots,
\al_{j-1}^{p^{m_j}})$.
\item[{\rm(2)}] Soit $\Lambda_j=\{(i_1,\dots, i_{j-1})$ tel que
$0\leq i_1<p^{m_1-m_j},\dots,0\leq i_{j-1}<p^{m_{j-1}-m_j}\}$,
alors $\{{(\al_1,\dots, \al_{j-1})}^{{p^{m_j}}\xi}$ tel que $\xi\in
\Lambda_j\}$ est une base de $k(K^{p^{m_j}})$ sur $k$.
\item[{\rm(3)}] Soient $n\in\mathbf{N}$ et $j$ le plus grand entier tel que
$m_j>n$. Alors $\{\al_{1}^{p^{n}},\dots, $ $\al_{j}^{p^{n}}\}$ est une
$r$-base canoniquement ordonn\'{e}e de $k(K^{p^n})/k$, et sa liste
des exposants est $(m_1-n,\dots, m_j-n)$}.
\end{itemize}
\end{pro}
\pre  cf. {\cite{Che-Fli2}}, p.
{140}, proposition {5.3}. \cqfd
\begin{pro} {\label{pr16}}
Soient $K_1/k$ et $K_2/k$ deux
sous-extensions purement ins\'ep\-a\-rables de $K/k$. $K_1$ et $K_2$
sont $k$-lin\'{e}airement disjointes si et seulement si
$o_j(K_1(K_2)/K_2)=o_j(K_1/k)$ pour tout $ j\in\mathbf{N}$.
\end{pro}
\pre cf. {\cite{Che-Fli}}, p.  {374},
 proposition  {7}. \cqfd
\begin{pro} {\label{pr17}} (Algorithme de la compl\'etion des r-bases) Soient K/k une
extension purement ins\'eparable finie, $G$ un $r$-g\'en\'erateur de $K/k$, et $\{\al_1,\dots, $ $\al_s\}$
un syst\`eme de $K$ tel que pour tout $j\in \{1,\dots, s\}$, $o(\al_j,$ $k($ $\al_1,\dots,\al_{j-1}))=o_j(K/k)$.
Pour toute suite $\al_{s+1}, \al_{s+2}, \dots,$  d'\'el\'ements de $G$  v\'erifiant
$o(\al_m,$ $k($ $\al_1,$ $\dots,\al_{m-1}))=\di\sup_{a\in G}(o(a,$ $k($ $\al_1,\dots,\al_{m-1})))$,
la suite $(\al_i)_{i\in \N^*}$ s'arr\^ete sur un plus grand entier $n$  tel que $o(\al_n, k(\al_1,\dots,\al_{n-1}))>0$. En particulier, $\{\al_1,\dots,$ $\al_n\}$ est une r-base canoniquement ordonn\'ee de $K/k$.
\end{pro}
\pre Cf. {\cite{Che-Fli2}}, p. {139},  proposition {1.3}.\cqfd
\sk

{\bf Cas o\`u $K/k$ est $q$-finie d'exposant non born\'e.}
Soit $K/k$ une extension $q$-finie. Rappelons que pour tout $n\in \N^*$, $k_n$ d\'esigne toujours $k^{p^{-n}}\cap K$. En vertu de la proposition \textcolor{blue}{\ref{pr14}},  pour tout $j\in \N^*$, la suite des entiers naturels $(o_j(k_n/k))_{n\geq 1}$ est croissante, et donc
$(o_j(k_n/k))_{n\geq 1}$ converge vers $+\infty$, ou $(o_j(k_n/k))_{n\geq 1}$ est stationnaire \`{a} partir d'un certain rang. Lorsque $(o_j(k_n/k))_{n\geq 1}$ est born\'{e}e,  par construction,  pour tout $t\geq j$, $(o_t(k_n/k))_{n\geq 1}$ est aussi born\'{e}e (et donc stationnaire).
\begin{df} Soient $K/k$ une extension $q$-finie  et $j$ un entier naturel non nul. On appelle le $j$-i\`{e}me exposant de $K/k$ l'invariant $o_j(K/k)=\di\lim_{n\rightarrow +\infty} (o_j(k_n$ $/k))$.
\end{df}
\begin{lem} {\label{lem2}} Soit $K/k$ une extension $q$-finie, alors  $o_s(K/k)$ est fini si et seul\-ement s'il existe un entier naturel $n$ tel que $di(k(K^{p^n})/k)<s$, et on a $o_s(K/k)=\inf\{m\in \N\,|\, di(k(K^{p^m})/k)<s\}$. En particulier, $o_s(K/k)$ est infini si et seulement si pour tout $m\in \N$, $di(k(K^{p^m})/k)\geq s$.
\end{lem}
\pre Pour simplifier l'\'{e}criture, on note $e_t=o_t(K/k)$ si $o_t(K/k)$ est fini. Compte tenu du  \cite{Che-Fli2}, p. {138}, lemme {1.3}, on v\'erifie aussit\^ot que $o_s(K/k)$ est infini si et seulement si pour tout $m\in \N$, $di(k(K^{p^m})/k)\geq s$, donc on se ram\`{e}ne au cas  o\`{u} $o_s(K/k)$ est fini. Par suite, il existe un entier $n_0$, pour tout $n\geq n_0$, $e_s=o_s(k_n/k)$. D'apr\`{e}s \cite{Che-Fli2} p. {138}, lemme {1.3},  $di(k({k_n}^{p^{e_s}})/k)<s$ et $di(k({k_n}^{p^{e_s-1}})/k)\geq s$. En vertu du th\'{e}or\`{e}me \textcolor{blue}{\ref{thm4}},  $di(k(K^{p^{e_s}})/k)<s$ et  $di(k({K}^{p^{e_s-1}})/k)\geq s$. Autrement dit,  $o_s(K/k)=\inf\{m\in \N\,|\, di(k(K^{p^m})/k)<s\}$. \cqfd \sk

Le r\'{e}sultat ci-dessous permet de ramener l'\'{e}tude des propri\'{e}t\'{e}s des exposants des extensions $q$-finies aux extensions finies par le biais  des cl\^{o}tures relativement parfaites.
\begin{thm} {\label{thm8}} Soit $K_r/k$ la cl\^{o}ture relativement parfaite de degr\'{e} d'irratio\-n\-alit\'{e} $s$ d'une extension $q$-finie $K/k$, alors on a :
\sk

\begin{itemize}{\it
\item[\rm{(i)}] Pour tout $t\leq s$, $o_t(K/k)=+\infty$.
\item[\rm{(ii)}] Pour tout $t> s$, $o_t(K/k)=o_{t-s}(K/K_r)$.}
\end{itemize}
\sk

\noindent En outre, $o_t(K/k)$ est fini si et seulement si $t> s$.
\end{thm}
\pre  Pour tout $t\in {\N}^*$, notons $e_t=o_t(K/K_r)$. Comme pour tout entier $e$, on a $k(K^{p^e})=K_r(K^{pe})=\di\bigcup_{n\in \N}k({k_n}^{p^e})$, donc  $s=di(K_r/k)\leq di(k(K^{p^e})/k)=di(k({k_n}^{p^e})/k)$ pour $n$ suffisament grand. D'apr\`{e}s le lemme \textcolor{blue}{\ref{lem2}}, on aura d'une part $o_t(K/k)=+\infty$ pour tout $t\leq s$, et
d'autrs part pour tout $n>s$, $di(K_r(K^{p^{e_{n-s}}})/k)=di(K_r/k)+di(K_r(K^{p^{e_{n-s}}})/K_r)<s+n-s=n$ et $di(K_r(K^{p^{e_{n-s}-1}})/k)=di(K_r/k)+di(K_r(K^{p^{e_{n-s}-1}})/K_r)\geq n$. Notamment, pour tout $n>s$, $o_n(K/k)=o_{n-s}(K/K_r)$. Toutefois, $o_n(K/k)$ est fini si et seulement si $n\leq s$.   \cqfd \sk

Voici une liste de cons\'equences imm\'ediates :
\begin{pro} {\label{pr18}} Soient $K$ et $L$ deux corps interm\'{e}diaires  d'une extension $q$-finie $M/k$. Pour tout  $j\in \N^* $, on a $o_j(L(K)/L)\leq o_j(K/k)$.
\end{pro}
\pre  Due au lemme \textcolor{blue}{\ref{lem2}},  et \`{a} l'in\'{e}galit\'{e} suivante r\'{e}sultant du corllaire \textcolor{blue}{\ref{cor4}}: $di(L(L^{p^{n}}, K^{p^n})/L)=di(L(K^{p^n})/L)\leq di (k(K^{p^n})/k)$ pour tout $n\in \N$. \cqfd
\begin{pro} {\label{pr19}} Etant donn\'ees des extensions $q$-finies $k\subseteq L\subseteq K$. Pour tout  $j\in \N^*$, on a $o_j(L/k)\leq o_j(K/k).$
\end{pro}
\pre Application imm\'ediate du lemme \textcolor{blue}{\ref{lem2}},  et de l'in\'{e}galit\'{e} suivante r\'{e}sultant du th\'{e}or\`{e}me \textcolor{blue}{\ref{thm4}} : $di(k(L^{p^{n}})/k) $ $\leq di (k(K^{p^n})/k)$ pour tout $n\in \N$. %\cqfd
\sk

Par ailleurs la taille d'une extension relativement parfaite reste invariant, \`a une extension finie pr\`es comme l'indique le r\'esultat suivant.
\begin{pro} {\label{pr20}}  Etant donn\'ee  une sous-extension $K/k$ relativement parfaite   d'une extension $q$-finie $M/k$. Pour toute sous-extension finie $L/k$ de $M/k$, on a $di(L(K)/L)=di(K/k)$.
\end{pro}
\pre  En vertu du corollaire \textcolor{blue}{\ref{cor4}}, il suffit de montrer que $di(L(K)/L)\geq di(K/k)$. Pour cela, on
pose d'abord $e=o_1(L/k)$ et $t=di(K/k)$. D'apr\`es le th\'eor\`eme \textcolor{blue}{\ref{thm8}}, pour tout $s\in \{1,\dots, t\}$, $o_s(K/k)=+\infty$, donc pour $n$ assez grand,  on aura $o_t(k_n/k)>e+1$, en outre $L\subseteq k_n$ et $di(k_n/k)=di(K/k)$. Soit $\{\al_1,\dots, \al_t\}$ une $r$-base canoniquement ordonn\'ee de $k_n/k$, s'il existe $s\in \{1,\dots, t\}$ tel que $\al_s\in L({k_n}^p)(\al_1,\dots, \al_{s-1})$, d'apr\`es la proposition \textcolor{blue}{\ref{pr14}}, on aura  $e<o_t(k_n/k)\leq o_s(k_n/k)=o(\al_s, k(\al_1,\dots,\al_{s-1}))\leq o_1(L({k_n}^p)(\al_1,\dots,\al_{s-1})/$ $k(\al_1,\dots,\al_{s-1}))\leq \di \sup(o_1(L/k), o_s(k_n/k)-1)=o_s(k_n/k)-1$, et donc $o_s(k_n/k)$ $\leq o_s(k_n/k)-1$, contradiction. D'o\`u, $\{\al_1,\dots, \al_t\}$ est une $r$-base de $L(k_n)/L$, et par suite, $t=di(K/k)=di(L(k_n)/L)\leq di(L(K)/L)$. \cqfd
%%%%%%%%%%%%%%%%%%%%%%%%%%%%%%%%%%%%%%%%%%%
\section{Extensions modulaires}

On rappelle qu'une extension
$K/k$ est dite modulaire si et seulement si pour tout
$n\in\mathbf{N}$, $K^{p^{n}}$ et $k$ sont $K^{p^{n}}\cap
k$-lin\'{e}airement disjointes. Cette notion a \'{e}t\'{e} d\'efinie
pour la premi\`{e}re fois par Swedleer dans {\cite{Swe}}, elle
caract\'{e}rise les extensions purement ins\'{e}parables, qui sont
produit tensoriel sur $k$ d'extensions simples sur $k$. Par ailleurs, toute $r$-base $B$ de $K/k$ telle que $K\simeq \otimes_k (\otimes_k k(a))_{a\in B}$ sera appel\'ee
$r$-base modulaire. En particulier, d'apr\`es le th\'eor\`eme de Swedleer,  si $K/k$ est d'exposant born\'e, il est \'equivalent de dire que :
\sk

\begin{itemize}{\it
\item[\rm{(i)}] $K/k$ admet une $r$-modulaire.
\item[\rm{(ii)}] $K/k$ est modulaire.}
\end{itemize}
\sk

Soient $m_j$ le $j$-i\`{e}me exposant d'une extension purement ins\'eparable finie $K/k$ et $\{\al_1 ,\dots,\al_n\}$
une $r$-base canoniquement ordonn\'{e}e de $K/k$,
donc d'apr\`es la proposition $\textcolor{blue}{\ref{pr15}}$, pour tout $j\in \{2,\dots,  n\}$, il existe des constantes uniques
$C_{\ep}\in k$ telles que ${\al_{j}}^{p^{m_j}}=\di\sum_{\ep \in \Lambda_j}{C_{\ep}}{(\al_1,
\dots, \al_{j-1})}^{p^{m_j}\ep}$, o\`u $\Lambda_j=\{(i_1,\dots, i_{j-1})$ tel que
$0\leq i_1<p^{m_1-m_j},\dots,0\leq i_{j-1}<p^{m_{j-1}-m_j}\}$.
Ces relations s'appellent les \'{e}quations de d\'{e}finition de $K/k$.
\sk

Le crit\`ere ci-dessous permet de tester la modularit\'e d'une extension.
\begin{thm} {\label{thm9}}
{\bf [Crit\`ere de modularit\'{e}]} Sous les notations ci-dessus,
les propri\'{e}t\'{e}s suivantes sont \'{e}quivalentes~:
\sk

\begin{itemize}{\it
\item[\rm{(1)}] $K/k$ est modulaire.
\item[\rm{(2)}] Pour toute $r$-base  canoniquement ordonn\'{e}e $\{\al_1,\dots , \al_n\}$
de $K/k$, les $C_{\ep}\in k \cap K^{p^{m_j}}$ pour tout $j\in \{2,\dots,  n\}$.
\item[\rm{(3)}] Il existe une $r$-base canoniquement ordonn\'{e}e  $\{\al_1,\dots , \al_n\}$
de $K/k$ telle que les $C_{\ep}\in k \cap K^{p^{m_j}}$ pour tout $j\in \{2,\dots,  n\}$.}
\end{itemize}
\end{thm}
\pre cf. {\cite{Che-Fli2}}, p.  {142}, proposition {1.4}. \cqfd
\begin{exe} Soient $Q$ un corps parfait de
caract\'eristique $p>0$, $k=Q(X,Y,$ $Z)$ le corps des fractions
rationnelles aux ind\'etermin\'ees $X,Y,Z$, et $K=k(\al_1,$ $\al_2)$ avec $\al_1=X^{p^{-2}}$ et
$\al_2=X^{p^{-2}}Y^{p^{-1}}+Z^{p^{-1}}$. On v\'erifie
aussit\^ot que
\sk

\begin{itemize}{\it
\item[$\bullet$] $o_1(K/k)=2$ et $o_2(K/k)=1$,
\item[$\bullet$] $\al_2^{p}=Y\al_1^p+Z$.}
\end{itemize}
\sk

Si $K/k$ est
modulaire, d'apr\`es le crit\`ere du modularit\'e, on aura $Y\in k\cap K^p$ et $Z\in k\cap K^p$, et donc $Y^{p^{-1}}$ et $Z^{p^{-1}}\in K$. D'o\`u
$k(X^{p^{-2}},Y^{p^{-1}}, Z^{p^{-1}})$ $\subset $ $K$, et par suite,
$di(k(X^{p^{-2}},Y^{p^{-1}},Z^{p^{-1}})/k)=3<$ $di(K/k)=2$, contradiction.
\end{exe}

Le r\'{e}sultat suivant est cons\'{e}quence imm\'{e}diate de la  modularit\'{e}.
\begin{pro} {\label{apr2}} Soient $m,n\in \Z$ avec $n\geq m$. Si
$K/k$ est modulaire, alors
 $K^{p^{m}}/k^{p^{n}}$  est modulaire.
\end{pro}

La condition $n\geq m$ assure $k^{p^{n}}\subset K^{p^{m}}$.
\begin{pro} {\label{proa1}} Soit $K/k$ une extension purement ins\'{e}parable finie (respectivement, et modulaire), et soit $L/k$ une
sous-extension de $K/k$ (respectivement, et modulaire) avec
$di (L/k) = s$. Si $K^p\subseteq L$,  il existe une r-base canoniquement ordonn\'{e}e (respectivement, et modulaire)
$ (\al_1, \al_2, \dots, \al_n)$ de $K/k$, et $e_1, e_2,\dots, e_s\in \{1,p\}$ tels que $({\al_1}^{e_1}, {\al_2}^{e_2},\dots, {\al_s}^{e_s})$ soit une r-base canoniquement
ordonn\'{e}e (respectivement, et modulaire) de $L/k$. De plus, pour tout $j\in \{1,\dots,  s\}$, on a $o_j (K/k) = o_j (L/k)$, auquel cas $e_j = 1$, ou $o_j (K/k) = o_j (L/k) + 1$, auquel cas
$e_j = p$.
\end{pro}
\pre Cf. {\cite{Che-Fli2}}, p. {146},
 proposition {8.4}. \cqfd \sk

Le th\'{e}or\`{e}me suivant de  Waterhouse joue un r\^{o}le important
dans l'\'{e}tude des extensions  modulaires (cf. \cite{Wat} Th\'{e}or\`{e}me 1.1).
\begin{thm} Soient $(K_{j})_{j\in I}$ une famille
de sous-corps d'un corps commutatif $\Omega$, et $K$ un autre sous-corps de $\Omega $. Si pour tout $j\in I$, $K$ et $K_j$ sont  $K\cap K_j$-lin\'{e}airement
disjoints,
alors $K$ et  $\di\bigcap_{j} K_{j}$ sont  $K\cap (\di\bigcap_{j} K_{j})$-lin\'{e}airement disjoint.
\end{thm}

Comme cons\'equence, la modularit\'e est stable par une intersection quelconque portant soit au dessus ou en dessous d'un corps commutatif. Plus pr\'ecis\'ement, on a :
\begin{cor} {\label{apr4}} Sous les m\^emes hypoth\`eses du th\'eor\`eme ci-dessus, on a :
\sk

\begin{itemize}{\it
\item[\rm{(i)}] Si pour tout $j\in I$, $K_{j}/k$ est modulaire, il en est de m\^eme de $\di\bigcap_{j} K_{j}/k$.
\item[\rm{(ii)}] Si pour tout $j\in I$, $K/K_j$ est modulaire, il en est de m\^eme de $K/\di\bigcap_{j} K_{j}$.}
\end{itemize}
\end{cor}

D'apr\`{e}s le th\'{e}or\`{e}me de Waterhouse, il existe une plus petite sous-extension $m/k$ de $K/k$ (respectivement une plus petite extension $M/K$) telle que $K/m$ (respectivement $M/k$) est modulaire. D\'esormais, on note $m=lm(K/k)$ et $M=um(K/k)$. Toutefois, l'extension $um(K/k)$ sera appel\'ee cl\^oture modulaire de $K/k$.

%%%%%%%%%%%%%%%%%%%%%commentaire
Comme application imm\'{e}diate de la proposition \textcolor{blue}{\ref{pr7}}, on a
\begin{pro} {\label{pr24}} Etant donn\'{e}es une $r$-base modulaire $B$ d'une extension modulaire $K/k$ et une famille $(e_a)_{a\in B}$ d'entiers tels que $0\leq e_a\leq o(a,k)$.  Soit $L=k(({a}^{p^{e_a}})_{a\in B})$, alors $L/k$ et $K/L$ sont modulaires, et $(B\setminus L)$, $(({a}^{p^{e_a}})_{a\in B}\setminus L)$ sont deux $r$-bases modulaires  respectivement de $K/L$ et $L/k$. En outre, pour tout $a\in B$, $o(a,L)=e_a$.
\end{pro}
\pre On se ram\`ene au cas  fini auquel le r\'esultat d\'ecoule de la transitivit\'e de la lin\'earit\'e disjointe.
En outre, pour toute partie $\{a_1,\dots, a_n\}$ d'\'el\'ement de $B$, $[L(a_1,\dots,a_n):L]=\di\prod_{i=1}^{n}p^{e_{a_i}}$. \cqfd \sk
%%%%%%%%%%%%%%%%%%%%%%%%%%%%%%%%%%%%%%%%%%%%%%%%%%%%%%%
%%%%%%%%%%%%%%%%%%%%%%%%%%%%%%%%%%%%%%%%%%

Dans la suite, pour tout $a\in B$, on pose $n_a=o(a,k)$. Consid\'erons maintenant les sous-ensembles $B_1$ et $B_2$ de $B$ d\'efinis par $B_1=\{a\in B\,|\, n_a>j\}$, $B_2=B\setminus B_1=\{a\in B\,|\, n_a\leq j\}$, ($j$ \'etant un entier ne d\'epassant pas $o(K/k)$).
\sk

Comme Application de la proposition pr\'ec\'edente, on a :
\begin{thm} {\label{thm11}} Sous les conditions pr\'ecis\'ees ci-dessus, pour tout entier $j< o(K/k)$, on a $k_j=k((a^{n_a-j})_{a\in B_1}, B_2)$.
\end{thm}
\pre   Comme $K/k$ est r\'eunion inductive d'extentions modulaires  engendr\'ees par des parties finies de $B$, et compte tenu de la distributivit\'e de l'intersection par rapport \`a la r\'eunion,   on peut supposer sans perdre de g\'en\'era\-l\-i\-t\'e que $K/k$ est finie d'exposant not\'e $e$. Soient $\{\al_1,\cdots,\al_n\}$ une $r$-base modulaire  et canoniquement ordonn\'{e}e de $K/k$, et $m_j$ le j-i\`{e}me exposant de $K/k$.  D\'{e}signons par   $s$ le plus grand entier tel que $m_s>j$, et $L=k(\al_{1}^{p^{m_1-j}},\dots,$ $\al_{s}^{p^{m_s-j}}, $ $\al_{s+1},\dots,$ $ \al_n)$. On v\'{e}rifie aussit\^{o}t que :
\sk

\begin{itemize}{\it
\item[\rm{(i)}] $L\subseteq k_j$,
 \item[\rm{(ii)}] $K\simeq k(\al_1)\otimes_k \dots\otimes_k k(\al_n)\simeq L(\al_1)\otimes_L \dots \otimes_L L(\al_s)$.}
\end{itemize}
\sk

\noindent Ainsi,  pour tout $x\in K$, il existe des constantes uniques  $C_{\ep}\in L$ telles que
$x=\di\sum_{\ep \in \Lambda}{C_{\ep}}{(\al_1,
\dots, \al_{s})}^{\ep}$, o\`u $\Lambda=\{(i_1,\dots, i_{s})$ tel que
$0\leq i_1<p^{m_1-j},\dots, 0\leq i_{s}<p^{m_{s}-j}\}$, et donc $x^{p^j}=\di\sum_{\ep \in \Lambda}{C_{\ep}}^{p^j}{({\al_1}^{p^j},
\dots, {\al_{s}}^{p^j})}^{\ep}$. Compte tenu de la proposition \textcolor{blue}{\ref{pr15}}, $x^{p^j}\in k$ (c'est-\`{a}-dire $x\in k_j$) si et seulement si $x^{p^j}={C_{(0,\dots, 0)}}^{p^j}$, ou encore $x=C_{(0,\dots,0)}$. Par suite $x\in k_j$ si et seulement si $x\in L$,  autrement dit  $k_j=L$. %\cqfd
\sk
%%%%%%%%%%%%%%%%%%%%%%%%%%%%%%%%%%%%

Comme cons\'equence imm\'ediate, dans le cas de modulaire le r\'esultat suivant exprime une propri\'et\'e de stabilit\'e de la taille d'un certains corps interm\'ediaires. Plus pr\'ecis\'ement,
\begin{cor}{\label{acor1}} Pour toute extension modulaire $K/k$, pour tout $n\in \N$, on a $di(k_n/k)=di(k_1/k)$. En particulier, $di(K/k)=di(k_1/k)$.
\end{cor}

Le r\'{e}sultat suivant est bien connu (cf. \cite{Kim}).
\begin{pro} {\label{pr23}} Soit $K/k$ une extension purement
ins\'{e}parable et modulaire~; soit pour tout $n\in\mathbf{N}$,
$K_n=k(K^{p^n})$. Alors $k_n/k$, $K/k_n$,  $K_n/k$ et $K/K_n$ sont modulaires.
\end{pro}
\begin{pro} {\label{apr1}} Soient $K_1$ et $K_2$ deux sous-extensions de $K/k$ telles que $K\simeq K_1\otimes K_2$. Si pour tout $i\in \{1,2\}$,  $K_i/k$ est modulaire,  il en est de m\^eme de $K/k$.
\end{pro}
\pre Cf. {\cite{Che-Fli4}}, p. {55},
 lemme {3.4}. \cqfd \sk

Le r\'{e}sultat suivant \'{e}tend trivialement les hypoth\`{e}ses de la proposition {3.3}, \cite{Mor-Vin},  p. {94}, ainsi que le th\'{e}or\`{e}me {3.2},  \cite{Dev1}, p. {289}.
Il utilise plus particuli\`{e}rement les propri\'{e}t\'{e}s du syst\`{e}me canoniquement g\'{e}n\'{e}rateur (pour plus d'information cf. \cite{Mor-Vin}, d\'efinition {1.32}, p. 29).
\begin{pro}  {\label{apr3}}Soient $K_1$ et $K_2$ deux sous-extensions de $K/k$ telles que $K\simeq K_1\otimes K_2$. Si $K/K_1$ est modulaire, et $K_2/k$ est d'exposant born\'e,  il existe une partie $B$ de $K$ telle que $K\simeq K_1\otimes_k (\otimes_k (k(\al)_{\al\in B})$.
\end{pro}
\pre D'abord, comme $K\simeq K_1\otimes_k K_2$, alors pour tout $i\in \N$, pour toute $r$-base $C$ de $k({K_2}^{p^i})/k$,  $C$ est aussi une $r$-base de $K_1({K_2}^{p^i})/K_1$.
Choisissons ensuite une $r$-base $B$ de $K_2/k$, comme $K_2/k$ est d'exposant fini,  alors $B$ est un $r$-g\'{e}n\'{e}rateur minimal de $K_2/k$. Soit $B_1,\dots, B_n$ une partition de $B$ v\'{e}rifiant $B_1=\{x\in B| \, o(x,k)=o_1(K_2/k)=e_1\}$ et, pour tout $1<i\leq n$,  $B_i=\{x\in B| \, o(x,k(B_1,\dots, B_{i-1}))=o_1(K_2/k(B_1,\dots,B_{i-1}))=e_i\}$. Il est clair que $e_1>\dots>e_n$, et  en vertu de la lin\'earit\'e disjointe,   pour tout $i\in \{1\dots,n\}$,  pour tout $x\in B_i$,  on a \'egalement $o(x,K_1(B_1,\dots, B_{i-1}))=o_1(K/K_1(B_1,\dots,B_{i-1}))\}=e_i$.  En particulier, pour tout $i\in \{2,\dots, n\}$, $(\di\prod_{\al}{(G)}^{{\al}p^{e_i}})_G$, o\`{u} $G$ est une partie finie d'\'{e}l\'{e}ments de $ B_1\cup\dots\cup B_{i-1}$ et les $\al$ sont convenablement choisis,  est une base respectivement de $k({K_2}^{p^{e_i}})$ sur $k$  et $K_1({K_2}^{p^{e_i}})=K_1({K}^{p^{e_i}})$ sur $K_1$. Notons $M_i$ cette base, et soit $x\in B_i$, il existe des $c_{\al}\in k$ uniques tels que $x=\di\sum_{\al}c_{\al}y_{\al}$, ($y_{\al}\in M_i$), en outre les $c_{\al}$ sont aussi uniques dans $K_1$.
 D'autre part, en vertu de la modularit\'e, pour tout $i\in \{1,\dots, n\}$,  $K^{p^{e_i}}$ et $K_1$ sont $K_1\cap K^{p^{e_i}}$-lin\'eairement disjointes. Comme $K_1({K_2}^{p^{e_i}})=K_1({K}^{p^{e_i}})$ et $M_i\subseteq K^{p^{e_i}}$, alors $M_i$ est aussi une base de $K^{p^{e_i}}$ sur $K_1\cap K^{p^{e_i}}$. En tenant compte de l'unicit\'e de l'\'ecriture de $x$ dans la base $M_i$,  on en d\'eduit par identification que les $c_{\al}\in k\cap K^{p^{e_i}}$, et donc ${B_i}^{p^{e_i}}\subseteq  k\cap K^{p^{e_i}}({K_1}^{p^{e_i}}({B_1}^{p^{e_i}},\dots, {B_{i-1}}^{p^{e_i}}))$ pour tout $i\in \{1\dots, n\}$. Par application du (\cite{Mor-Vin}, proposition {3.3}, p. {94}), il existe une sous-extension modulaire $J/k$ d'exposant fini de $K/k$ telle que $K\simeq K_1\otimes_k J$. Ainsi, le r\'esultat d\'ecoule imm\'ediatement du th\'eor\`eme de Swedleer.  \cqfd \sk
%%%%%%%%%%%%%%%%%%%%%%%%%%%%%%%%%%%

Dans le cas fini, le r\'esultat suivant g\'en\'eralise la proposion ci-dessus.
\begin{pro} {\label{ajpr1}} Soient $K_1$ et $K_2$ deux corps interm\'{e}diaires ; $k$-lin\'{e}airement disjoints d'une extension purement ins\'eparable finie $L/k$ avec $di(L/K_1)=di(K_2$ $/k)=n$.
Soit $s$ le plus petit entier tel que $o_s(K_2/k)=o_n(K_2/k)$. Si $L/K_1$ est modulaire, il existe une $r$-base $\{\al_1,\dots, \al_n\}$  canoniquement ordonn\'{e}e de
$K_1(K_2)/K_1$ v\'{e}rifiant $K_1(K_2)\simeq K_1\otimes k(\al_1,\dots, \al_s)\otimes_k k(\al_{s+1})\otimes_k \dots \otimes_k k(\al_n)$.
\end{pro}
\pre Pour simplifier l'\'ecriture, pour tout $j\in\{1,\dots, n\}$, on note  $o_j (K_2/k)$ $ = e_j$, et $K=K_1(K_2)$ . Soit $\{\al_1,\dots,\al_n\}$ une $r$-base canoniquement ordonn\'{e}e de $K_2/k$.
Compte tenu de la proposition \textcolor{blue}{\ref{pr16}}, $\{\al_1,\dots,\al_n\}$ est aussi une $r$-base canoniquement ordonn\'{e}e de $K/K_1$, et pour tout $j\in\{1,\dots, n\}$,   $o_j (K/K_1) = e_j$.
D'apr\`{e}s la proposition \textcolor{blue}{\ref{pr15}}, pour tout $i\in \{s,\dots, n\}$, il existe des constantes uniques $C_{\ep}^{i}\in k$ telles que  ${\al_i}^{p^{e_n}}=\di\sum_{\ep\in \Lambda_{s-1}} C_{\ep}^{i}{(\al_1\dots\al_{s-1})}^{p^{\ep}}$ ($*$). En vertu de la proposition \textcolor{blue}{\ref{pr16}}, pour tout $i\in \{s\dots,n\}$,
l'\'{e}quation de d\'{e}finition de $\al_i$ par rapport \`{a} $K_1(\al_1,\dots,\al_{s-1})$ est aussi d\'efinie par la relation ($*$) ci-dessus.
Comme $L/K_1$ est modulaire, en se servant du crit\`{e}re de modularit\'{e}, pour tout $(i,\ep)\in \{s,\dots, n\}\times \Lambda_{s-1}$,
on aura ${(C_{\ep}^{i})}^{p^{-e_n}}\in L$.    Posons ensuite, $F=k({(C_{\ep}^{i})}^{p^{-e_n}})$ o\`u $(i,\ep)$ parcourt l'ensemble $\{s,\dots, n\}\times \Lambda_{s-1}$,
et $H=K_1(F)(\al_1,\dots,\al_{s-1})$. Il est clair que $o_1(F/k)\leq e_n$, et $K\subseteq H\subseteq L$. De plus, d'apr\`{e}s le th\'eor\`eme \textcolor{blue}{\ref{thm4}} et la proposition \textcolor{blue}{\ref{pr14}},
$n=di(K/K_1)\leq di(H/K_1)\leq di(L/K_1)=n$, et pour tout $i\in \{s, \dots, n\}$, $e_n=o_i(K/K_1)\leq o_i(H/K_1) \leq e_n$. Il en r\'{e}sulte que  $di(H/K_1)=n$,
et pour tout $i\in \{s, \dots, n\}$, $e_n=o_i(H/K_1)$. Comme $e_{s-1}>e_s=e_n$, d'apr\`es l'algorithme de la compl\'{e}tion des $r$-bases, il existe des \'el\'ements $b_s,\dots, b_n \in  F$
tels que $\{\al_1\dots,\al_{s-1},b_s,\dots,b_n\}$ soit une $r$-base canoniquement ordonn\'{e}e de $H/K_1$. En particulier, on aura :
\sk

\begin{itemize}{\it
\item[$\bullet$] Pour tout $i\in \{1 ,\dots, s-1\}$, $e_i=o_i(H/K_1)=o_i(K_1(\al_1,\dots,\al_{s-1})/K_1)= o_i(k(\al_1,\dots,\al_{s-1})/k)$.
\item[$\bullet$] Pour tout $j\in \{s,\dots,n\}$,  $e_n=o_j(H/K_1)=o(b_j, K_1(\al_1\dots,
\al_{s-1}, b_s,\dots, $ $b_{j-1})) \leq o(b_j, k(b_s,\dots,b_{j-1})/k)\leq o_1(F/k)\leq e_n$,
et donc $e_n=o_j(H/K_1)$ $=o_j(k(b_s,\dots,b_n)/k)$.}
\end{itemize}
\sk

\noindent D'o\`u, $H=K\simeq K_1\otimes k(\al_1,\dots, \al_{s-1})\otimes_k k(b_s)\otimes_k\dots\otimes_k k(b_n)$. \cqfd
%%%%%%%%%%%%%%%%%%%%%%%%%%%%%%%%µµµµµµµµµµµµµµµµµµµµµµ
%%%%%%%%%%%%%%%%%%%%%%%%%%%%%%%%µµµµµµµµµµµµµµµµµµµµµµ
\section{Extensions \'equiexponentielles }
\begin{pro} {\label{pr26}} Soit $K/k$ une extension purement ins\'{e}parable d'exposant $e$. Les  assertions suivantes sont \'equivalentes :
\sk

\begin{itemize}{\it
\item[\rm{(1)}] Il existe une $r$-base $G$ de $K/k$ v\'{e}rifiant $K\simeq \otimes_k (k(a))_{a\in G}$, et pour tout $a\in G$, $o(a,k)=e$.
\item[\rm{(2)}] Toute  $r$-base $G$ de $K/k$ satisfait $K\simeq \otimes_k (k(a))_{a\in G}$, et $o_1(K/k)=e$.
\item[\rm{(3)}] Il existe une $r$-base $G$ de $K/k$ telle que pour tout $a\in G$, $o(a,k(G\setminus \{a\}))=o(a,k)=e$.
\item[\rm{(4)}] Pour toute $r$-base $G$ de $K/k$, pour tout $a\in G$, $o(a,k(G\setminus \{a\}))=o(a,k)=e$.}
\end{itemize}
\end{pro}
\pre D'apr\`es le th\'eor\`eme de la $r$-base incompl\`ete, on se ram\`{e}ne au cas o\`{u} $K/k$ est finie auquel cas $[K:k]=p^{en}$, o\`u $e=o_1(K/k)$ et $n=di(K/k)$, et en vertu de la proposition \textcolor{blue}{\ref{pr16}}, le r\'esultat est imm\'ediat. \cqfd
\begin{df} Une extension qui v\'{e}rifie l'une des conditions de la proposition ci-dessus est dite \'equiexponentielle d'exposant $e$.
\end{df}

Il est clair  que toute extension \'equiexponentielle est modulaire. De plus, on v\'erifie aussit\^ot qu'il est \'{e}quivalent de dire que :
\sk

\begin{itemize}{\it
\item[\rm{(1)}] $K/k$ est \'equiexponentielle d'exposant $e$.
\item[\rm{(2)}] Il existe une $r$-base $G$ de $K/k$, pour toute partie finie $G_1$ de $G$, on a $k(G_1)/k$ est \'equiexponentielle d'exposant $e$.
\item[\rm{(3)}] Pour toute $r$-base $G$ de $K/k$, pour toute partie finie $G_1$ de $G$, on a $k(G_1)/k$ est \'equiexponentielle d'exposant $e$.}
\end{itemize}
\begin{pro} {\label{thm13}}
Pour toute extension $K/k$ relativement parfaite et modulaire, pour tout entier $n$, $k_n /k$ est \'equiexponentielle d'exposant $n$.
\end{pro}
\pre D'apr\`es le th\'eor\`eme \textcolor{blue}{\ref{thm11}}, il suffit de montrer que $k({k_n}^p)=k_{n-1}$. Compte tenu de la modularit\'e de $K/k$,  $K^{p^n}$ et $k$ sont $k\cap K^{p^n}$-lin\'{e}airement disjointes pour tout $n\geq 1$, et en vertu de la transitivit\'{e} de la lin\'{e}arit\'{e} disjointe, $k^{p^{n-1}}(K^{p^n})$ et $k$ sont $k^{p^{n-1}}(k\cap K^{p^{n}})$-lin\'{e}airement disjointes. Or $K/k$ est relativement parfaite, donc $k^{p^{n-1}}(K^{p^n})=K^{p^{n-1}}$, et par suite $k\cap K^{p^{n-1}}=k^{p^{n-1}}(k\cap K^{p^n})$, ou encore $k({k_n}^p)=k_{n-1}$. \cqfd \sk

Le r\'{e}sultat suivant, rapporte plus de pr\'{e}cision \`a la proposition \textcolor{blue}{\ref{thm13}} dans le cas des extensions $q$-finies, notamment aux extensions finies.
\begin{pro} {\label{pr27}}
Soit $K/k$ une extension purement ins\'{e}parable de degr\'e d'irr\-a\-t\-ionalit\'e $t$, relativement parfaite et modulaire (respectivement finie
et \'{e}quiexp\-o\-n\-entielle). Soient $n$ et $m$ deux entiers naturels tels que $n< m\hbox{ }(\hbox{respectivement, } n$ $< o_1(K/k))$. Les
propri\'{e}t\'{e}s suivantes sont v\'{e}rifi\'{e}es:
\sk

\begin{itemize}{\it
\item[{\rm (1)}] $di(k_{m}/k_{n})=t$.
\item[{\rm (2)}] $k_{m}/k_{n}$ est \'{e}quiexponentielle d'exposant
$m-n$;
\item[{\rm (3)}] $k_{n}^{p^{-(m-n)}}\cap K=k_{m}\;$ et
$\;k(k_{m}^{p^{m-n}})=k_{n}$.}
\end{itemize}
\sk

\noindent En particulier, pour tout $n\in \N$,  on a $[k_{n},k]=p^{nt}$.
\end{pro}
\pre cf. \cite{Che-Fli2}, p. {147}, proposition {9.4}. \cqfd \sk

Comme cons\'{e}quence imm\'{e}diate, on a :
\begin{cor} {\label{cor10}}
Si $K/k$ est une extension \'equiexponentielle d'exposant $e$, alors:
\sk

\begin{itemize}{\it
\item[\rm{(i)}] Pour tout $i\in \{1,\dots, e\}$, $k_i/k$ et $K/k_i$ sont \'equiexponentielles  d'exposant respectivement $i$ et $e-i$.
\item[\rm{(ii)}] Pour tout $i\in\{ 1,\dots, e\}$, $k(K^{p^{i}})/k$ et $K/k(K^{p^i})$ sont \'equiexponentielles  d'exposant respectivement $e-i$ et $i$.}
\end{itemize}
\end{cor}
\pre Imm\'ediat. \cqfd \sk

Le th\'eor\`eme ci-dessus reproduit dans un cadre plus \'etendu le corollaire {4.5} qui se trouve dans \cite{Dev1}, p. {292}, et pour plus d'information au sujet d'extraction des $r$-bases modulaires, on se r\'{e}f\`{e}re aux \cite{Dev1} et \cite{Dev-Mor2}.
\begin{thm} {\label{thm12}} Soient $k\subseteq L\subseteq K$ des extensions purement ins\'{e}parables telles que $K/k$ est \'equiexponentielle d'exposant $e$. Si $K/L$ est modulaire, il existe une $r$-base $G$ de $K/k$ telle que $\{a^{p^{o(a,L)}}| \, {a\in G}$ et $o(a,L)<e\}$ est une $r$-base modulaire de $L/k$.
\end{thm}
\pre Comme $K/L$ est modulaire d'exposant fini,  il existe une $r$-base $B_1$ de $K/L$ telle que $K\simeq \otimes_L (\otimes_L L(a))_{a\in B_1})$, (*). Pour des raisons d'\'{e}criture, pour tout $a\in B_1$, on pose $e_a=o(a,L)$ et $C=(a^{p^{e_a}})_{a\in B_1}$. Soit $B_2$ une partie de $L$ telle que $B_2$ est une $r$-base de $L(K^p)/k(K^p)$. Compte tenu de la transitivit\'e de $r$-ind\'ependance, $B_1\cup B_2$ est aussi une $r$-base de $K/k$. Dans la suite, notons $M=k(C,B_2)$. Il est clair que $M\subseteq L$, de plus, comme $K/k$ est \'equiexponentielle, on aura $K\simeq \otimes_k (\otimes_k k(a))_{a\in B_1\cup B_2}$. En vertu de la transitivit\'{e} de la lin\'{e}arit\'{e} disjointe, $K\simeq \otimes_M (\otimes_M M(a))_{a\in B_1}$, (**). En particulier,  d'apr\`{e}s les relations (*) et (**),  pour toute famille finie $\{a_1,\cdots, a_n\}$ d'\'el\'ements de $B_1$,  $L(a_1,\dots, a_n)\simeq L(a_1)\otimes_L\dots \otimes_L L(a_n)$ et $M(a_1,\dots, a_n)\simeq M(a_1)\otimes_M\dots\otimes_M M(a_n)$. Par application de la proposition \textcolor{blue}{\ref{pr16}}, on a successivement $ [L(a_1,\dots, a_n) :L]=\di\prod_{i=1}^{n}p^{e_{a_i}}$ et $[M(a_1,\dots, a_n) :M]=\di\prod_{i=1}^{n}p^{e_{a_i}}$, ou encore $L$ et $K$ sont $M$-lin\'{e}airement disjointes. D'o\`{u} $L=L\cap K=M$. \cqfd
 %%%%%%%%%%%%%%%%%%%%%%%%%%%%%%%%%%%%%%%%%%%%%%%%%%%
%%%%%%%%%%%%%%%%%%µµµµµµµµµµµµµµµµµµµµµµµµµµµµµµµµµ
\section{$q$-finitude et modularit\'{e}}

Soit $K/k$ une extension $q$-finie d'exposant non born\'e. Dans tout ce qui suit, nous utilisons les notations suivantes : $k_j= k^{p^{-j}}\cap K$, $U_{s}^j(K/k)=j-o_s(k_j/k)$, et  $Ilqm(K/k)$ d\'{e}signe
le premier entier $i_0$ pour lequel la suite  $(U_{i_0}^{j}(K/k))_{j\in\mathbf{N}}$ est non born\'{e}e.
Le r\'{e}sultat ci-dessus est une application imm\'{e}diate de la proposition \textcolor{blue}{\ref{pr14}}.
\begin{pro} {\label{pr25}} Etant donn\'ee une extension $q$-finie $K/k$.
 Pour tout entier $s$,  la suite $(U_{s}^j(K/k))_{j\in\mathbf{N}}$ est croissante.
\end{pro}
\pre
Comme ${k_{n+1}}^p\subseteq k_n$, il est clair que $o_s(k_n/k)\leq o_s(k_{n+1}/k) \leq o_s(k_n/k)+1$,  et donc $n+1-
o_s(k_{n+1}/k)\geq n-o_s(K_n/k)$~; c'est-\`{a}-dire  la suite
$(U_{s}^j(K/k))_{j\in\mathbf{N}}$ est croissante. \cqfd \sk

En outre, on v\'erifie aussit\^ot que :
\sk

\begin{itemize}{\it
\item[\rm{(i)}] Pour tout $s\geq Ilqm(K/k)$, $\di\lim_{n \rightarrow
+\infty}(U_{s}^n (K/k))=+\infty$.
\item[\rm{(ii)}]  Pour tout $s<Ilqm(K/k)$,  la suite
$(U_{s}^j(K/k))_{j\in\mathbf{N}}$ est born\'{e}e ;}
\end{itemize}

\noindent et par suite, pour tout
$n\geq \di\sup_{j\in\mathbf{N}}(\di\sup (U_{s}^{j}(K/$ $k$ $)))_{s<Ilqm(K/k)}$, on
a $U_{s}^n(K/k)=U_{s}^{n+1}(K/k)$. Autrement dit,
$o_s(k_{n+1}/k)=o_s(k_n/k)+1$.
\sk

Dans toute la suite,  on pose $e(K/k)=\di\sup_{j\in\mathbf{N}}(\di\sup (U_{s}^{j}(K/k)))_{s<Ilqm(K/k)}$, et pour tout $(s,j)\in {\N}^*\times {\N}^{*}$, $e_{s}^{j}=o_s(k_j/k)$
\begin{thm} {\label{thm10}}  Soit $K/k$ une extension $q$-finie, avec $t=di(rp(K/k)/k)$. Les affirmations suivantes sont \'equivalentes ;
\sk

\begin{itemize}{\it
\item[\rm{(1)}] $K/k$ est modulaire sur une extension finie de $k$.
\item[\rm{(2)}]  Pour tout $s\in\{1,2,\dots,t\}$,  la suite  $(U_{s}^{j}(K/k))_{j\in\mathbf{N}}$  est born\'{e}e.
\item[\rm{(3)}] $Ilqm(K/k)=t+1$.}
\end{itemize}
\end{thm}
\pre
Il est clair que $(2) \Leftrightarrow (3)$. Par ailleurs, compte tenu de la proposition \textcolor{blue}{\ref{pr10}}, il existe un entier $j_0$ tel que $K/k_{j_0}$ est relativement parfaite et $k_{j_0}(rp(K/k))=K$, et d'apr\`es la proposition \textcolor{blue}{\ref{pr20}}, on aura $di(K/k_{j_0})= di(rp(K/$ $k)/k)=t$. Supposons ensuite que la condition $(1)$ est v\'{e}rifi\'{e}e. On distingue deux cas :
\sk

 Si $K/k$ est modulaire, en vertu de la proposition \textcolor{blue}{\ref{pr27}}, pour tout $j\geq j_0$, on a $k_j/k_{j_0}$ est  \'{e}quiexponentielle d'exposant $j-j_0$ et $di(k_j/k_{j_0})=t$. D'o\`{u} pour tout $s\in \{1,\dots,t\}$, on a $U_{s}^{j}(K/k)=U_{s}^{j+1}(K/k)$.
 \sk

Si $K$ est modulaire sur une extension finie $L$ de $k$,  compte tenu de la finitude de  $L/k$, il existe un entier naturel $n$ tel que $L\subseteq k_n$. Par suite, $L^{p^{-j}}\cap K\subseteq k_{n+j}$, et donc $U_{s}^{n+j}(K/k)\leq n+U_{s}^{j}(K/L)$. D'o\`{u}, la suite $(U_{s}^{j}(K/k))_j$ est stationnaire pour tout $s\in\{1,\dots, t\}$.
\sk

 Inversement, si la  condition $(2)$ est v\'{e}rifi\'{e}e, il existe $m_0\geq \sup(e(K/k),j_0)$,  pour tout $j\geq m_0$, pour tout $s\in \{1,\dots,t\}$,  on a $o_s(k_{j+1}/k)=o_s(k_j/k)+1$ (et $di(k_j/k_{m_0})=t$). Par suite, $k_j/k_{j_0}$ est \'{e}quiexponentielle, donc modulaire. D'o\`{u} $K=\di\bigcup_{j>m_0} k_j$ est modulaire sur $k_{j_0}$. \cqfd
 \begin{thm} {\label{thm27}}
La plus petite sous-extension $M/k$ d'une extension $q$-finie $K/k$ telle que
$K/M$ est modulaire n'est pas triviale ($M\not= K$). Plus
pr\'{e}cis\'{e}ment, si $K/k$ est d'exposant non born\'e, il en est de m\^eme de $K/M$.
\end{thm}
%%%%%%%%%%%%%%%%%%%%%%%%%%%%%%%%%%%%%%%%%%%%%%%%%%%
\pre
Le cas o\`{u} $K/k$ n'est pas relativement parfaite (en particulier le cas fini)  est trivialement \'{e}vident, puisque $K/k(K^p)$ est modulaire. Ainsi, on est amen\'{e} \`{a} consid\'erer  que $K/k$ est relativement parfaite d'exposant non born\'{e}. On  emploiera ensuite un raisonnement par r\'{e}currence sur $di(K/k)=t$. Si $t = 1$, ou encore si $K/k$ est $q$-simple, il est imm\'{e}diat que  $K/k$ est modulaire.
Supposons maintenant que $t > 1$,  si  $Ilqm(K/k)=t+1$, en vertu du th\'eor\`eme \textcolor{blue}{\ref{thm10}}, $M/k$ est finie, et donc $K/M$ est d'exposant non born\'e. Si  $Ilqm(K/k)\leq t$, pour tout $j>e(K/k)$, pour tout $s\in[1; i-1]$ o\`u $i=Ilqm(K/k)$, on a $e^{j+1}_{s}=e^{j}_{s}+ 1$. Comme $k^{p}_{j+1}\subseteq  k_j$, d'apr\`es la proposition \textcolor{blue}{\ref{proa1}}, il
existe une r-base canoniquement ordonn\'ee $(\al_1,\dots, \al_n)$ de $k_{j+1}/k$, il existe $\ep_i,\dots,\ep_t\in \{1,p\}$  tels que $(\al_{1}^{p},\dots,\al_{i-1}^{p},\al_{i}^{\ep_i} \dots,\al_{t}^{\ep_t})$  est une r-base canoniquement ordonn\'{e}e de $k_j/k$. Dans la suite, pour tout $j>e(K/k)$, notons $K_j = k(k_{j}^{p^{e^{j}_{i}}})$. D'une part,  $K_j=k(\al_{1}^{p^{e^{j}_{i}+1}},\dots,\al_{i-1}^{p^{e^{j}_{i}+1}})$ et $K_{j+1}= k(\al_{1}^{p^{e^{j+1}_{i}}},\dots,\al_{i-1}^{p^{e^{j+1}_{i}}})$. D'autre part, on a $e^{j+1}_{i}=e^{j}_{i}+\ep $, avec $\ep=0$ ou $1$, cela conduit \`a $K_j\subseteq K_{j+1}$. Toutefois,  par d\'{e}finition de $Ilqm(K/k)$, on a $1 + e^{j}_{i }> e^{j+1}_{i}$ (c'est-\`a-dire $e^{j}_{i} = e^{j+1}_{i}$) pour une infinit\'{e} de valeurs de $j$. Pour ces valeurs, on a $di(K_{j+1}/k) = i -1$, sinon d'apr\`es le lemme \textcolor{blue}{\ref{lem2}},  $e^{j+1}_{i} = e^{j+1}_{i-1} = 1 + e^{j}_{i-1} = e^{j}_{i}$ , et donc $e^{j}_{i} > e^{j}_{i-1}$, ce qui contredit la d\'efinition des exposants.
Comme $(di(K_j/k))_{j>e(K/k)}$ est une suite croissante d'entiers born\'{e}e par $di(K/k)$, donc elle stationne sur $Ilqm(K/k)-1$. De plus, $K_j \not= K_{j+1}$, en effet si $K_j =K_{j+1} =k( K^{p}_{j+1})$,  comme $K_{j+1}/k$ est d'exposant born\'e, on aura $K_{j+1} = k$, ce qui est absurde.
Posons ensuite $H =\di\bigcup_{j>e(K/k)} K_j$. On v\'erifie aussit\^ot que $H/k$ est d'exposant non born\'{e} et $di(H/k) = i -1$, de plus $H/k$ est relativement parfaite car $k(K^{p}_{j+1})= K_j$ pour une infinit\'{e} de $j$. Par ailleurs, d'apr\`{e}s les corollaires \textcolor{blue}{\ref{cor8}} et \textcolor{blue}{\ref{ccor1}}, $di(K/H)<t$ et $K/H$ est d'exposant non born\'e.
Compte tenu de l'hypoth\`{e}se de r\'{e}currence appliqu\'ee \`a $K/H$, on aura $K$ est modulaire sur une extension $M'$ de $H$ avec
$K/M'$ est d'exposant non born\'e; comme $M \subseteq M'$, alors $K/M$ est aussi d'exposant non born\'e. \cqfd \sk

Une version \'equivalente de ce  r\'{e}sultat  se trouve dans \cite{Che-Fli2}. Toutefois, le th\'eor\`eme ci-dessus peut tomber en d\'efaut lorsque l'hypoth\`ese de la $q$-finitude n'est pas v\'erifi\'ee comme le montre le contre-exemple ci-dessus
\begin{exe}
Soient $Q$ un corps parfait de caract\'{e}ristique $p>0$, et $(X,(Y_i$ $)_{i\in\N^*}, $ $(Z_i)_{i\in \N^*}, (S_i)_{i\in\N^*})$ une famille alg\'{e}briquement ind\'{e}pendante sur $Q$. Soit $k=Q(X,(Y_i)_{i\in\N^*}, (Z_i)_{i\in \N^*}, (S_i)_{i\in\N^*})$ le corps des fractions rationnelles aux ind\'{e}\-t\-ermin\'{e}es $(X,(Y_i)_{i\in\N^*}, (Z_i)_{i\in \N^*}, (S_i)_{i\in\N^*})$.
Posons ensuite :
\sk

\begin{itemize}{\it
\item[] $K_1=\di\bigcup_{n\geq 1} k(\te_{1,n})$, avec $\te_{1,1}={X}^{p^{-1}}$ et $\te_{1,n}={\te_{1,n-1}}^{p^{-1}}$ pour tout entier $n>1$.
\item[] $K_2=\di\bigcup_{n\geq 1} K_1(\te_{2,n})$, o\`u $\te_{2,1}={Z_1}^{p^{-1}}\te_{1,2}+{Z_2}^{p^{-1}}$, et  pour tout $n>1$, $\te_{2,n}={Z_1}^{p^{-1}}\te_{1,2n}+{\te_{2,n-1}}^{p^{-1}}$.
\item[]  Par r\'{e}currence, on pose $K_j=\di\bigcup_{n\geq 1} K_{j-1}(\te_{j,n})$, o\`u $\te_{j,1}={Z_{j-1}}^{p^{-1}}\te_{j-1,2}$ $+{Z_j}^{p^{-1}}$, et pour tout $n>1$, $\te_{j,n}={Z_{j-1}}^{p^{-1}}\te_{j-1,2n}$ $+{\te_{j,n-1}}^{p^{-1}}$.}
\end{itemize}
\sk

Enfin, on note  $K=\di\bigcup_{j\in\N^*} K_j$, et  par conventient on pose $K_0=k$, et pout tout $i\in \N$, $\te_{i,0}=0$.
Comme pour tout $j\in \N^*$, $K_j\subseteq K_{j+1}$, alors $K$ est un corps commutatif.
\end{exe}
\begin{thm} {\label{thm14}} Sous les conditions ci-dessus, la plus petite sous-extension $m$ telle que $K/m$ est modulaire est triviale, c'est-\`a-dire $lm(K/k)=K$
\end{thm}

Pour la preuve de ce th\'{e}or\`{e}me, on se servira des r\'{e}sultats  suivants :
\begin{lem} {\label{lem3}} Sous les m\^{e}mes conditions ci-dessus, pour tout $(j,n)\in \N\times \N^*$, $K_j(\te_{j+1,n})=K_j({\te_{j+1,{n+1}}}^{p})$ et $\te_{j+1,1}\not\in K_j$. En particulier,  $o(\te_{j+1,n},K_j)=n$.
\end{lem}
\pre Il est trivialement \'evident que $K_j(\te_{j+1,n})=K_j({\te_{j+1,{n+1}}}^{p})$ pour tout $(j,n)\in \N\times \N^*$.  Pour achever la preuve, il suffit de remarquer que $K_j\subseteq k({X}^{p^{-\infty}},$ ${Z_1}^{p^{-\infty}},$ $\dots, {Z_j}^{p^{-\infty}})$ et $k(\te_{j+1,n},{X}^{p^{-\infty}},{Z_1}^{p^{-\infty}},$ $\dots, {Z_j}^{p^{-\infty}})= k({Z_{j+1}}^{p^{-n}}, {X}^{p^{-\infty}},$ ${Z_1}^{p^{-\infty}},$ $\dots, {Z_j}^{p^{-\infty}})$,  et donc, pour tout $n\in \N^*$, $n=o(\te_{j+1,n},k({X}^{p^{-\infty}},{Z_1}^{p^{-\infty}},$ $\dots, {Z_j}^{p^{-\infty}}$ $))\leq o(\te_{j+1,n},K_j)\leq n $. \cqfd \sk

Comme cons\'equence imm\'ediate, pour tout $j\in \N^*$, $K_j/K_{j-1}$  est $q$-simple d'exposant non born\'e. En particulier, $di(K_j/k)=j$.
\begin{lem} {\label{lem4}} Pour tout $i\in \N^*$, la famille $(Z_i,(S_j)_{j\in \N^*})$ est $r$-libre sur $K^p$.
\end{lem}
\pre
Puisque pour tout $i\geq 1$, ${S_i}^{p^{-1}}\not\in k({X}^{p^{-\infty}}, ({Z_{j}}^{p^{-\infty}})_{j\geq 1})({S_1}^{p^{-1}},\dots, $ ${S_{i-1}}^{p^{-1}})=K({Z_1}^{p^{-\infty}})({S_1}^{p^{-1}},\dots, {S_{i-1}}^{p^{-1}})$, il suffit de montrer que ${Z_i}\not\in K^p$; ou encore ${Z_i}^{p^{-1}}\not\in K$. Par construction, pour tout $j\in \{1,\dots, n\}$, on a $\te_{j,1}={Z_{j-1}}^{p^{-1}}\te_{j-1,2}$ $+{Z_j}^{p^{-1}}$ avec $K_n$ contient $\te_{j,1}$ et $\te_{j-1,2}$, et donc s'il existe $n>i$ tel que ${Z_i}^{p^{-1}}\in K_n$, par it\'eration, on aura  ${Z_{i-1}}^{p^{-1}},\dots, {Z_1}^{p^{-1}}\in K_n$ et ${Z_{i+1}}^{p^{-1}},\dots, {Z_n}^{p^{-1}}\in K_n$. Par suite, d'apr\`es le th\'eor\`eme \textcolor{blue}{\ref{thm4}}, $di(k(X^{p^{-1}}, {Z_1}^{p^{-1}},$ $\dots, {Z_n}^{p^{-1}})/k)\leq di(K_n/k)$, ou encore $n+1\leq n$, absurde.
D'o\`u pour tout $n\in \N^*$, ${Z_i}^{p^{-1}}\not\in K_n$, et comme $K$ est r\'eunion  de la famille croissante d'extensions $(K_n)_{n\in \N^*}$, alors  ${Z_i}^{p^{-1}}\not\in K$. \cqfd \sk

\noindent{}{\bf Preuve du th\'{e}or\`{e}me \textcolor{blue}{\ref{thm14}}.} Posons $m=lm(K/k)$. En utilisant un raisonnement par r\'ecurrence on va montrer que  $K_i\subseteq  m$  pour tout $i\in {\N}$, et par suite obtenir $K=m$. Il est imm\'ediat que  $K_0=k\subseteq m$, donc le r\'esultat est v\'erifi\'e pour le rang $0$. Soit $i\in \N^*$, supposons par application de  l'hypoth\`ese de r\'ecurrrence que $K_i\subseteq m$. S'il existe un entier naturel $s$ tel que $\te_{{i+1},s}\not \in m$,  d\'esignons par $n$ le plus grant entier tel que $\te_{{i+1},n}\in m$. D'o\`{u} pour tout $t\in \{0, \dots, n\}$, $\te_{{i+1},t}\in m$ et $\te_{{i+1},n+1}\not\in m$, en outre $\te_{{i+1},2n}^{p^n}\in m$, et ${\te_{{i+1},2(n+1)}}^{p^{n+1}}$ $\not\in m$. Il en r\'{e}sulte que le syst\`{e}me $({\te_{{i+1},2(n+1)}}^{p^{n+1}},1)$ est libre sur $m$, en particulier, il en est de m\^{e}me sur $m\cap K^{p^{n+1}}$. Compl\'{e}tons ce syst\`{e}me en une base $B$ de $K^{p^{n+1}}$ sur $m\cap K^{p^{n+1}}$. Comme $K^{p^{n+1}}$ et $m$ sont $m\cap K^{p^{n+1}}$-lin\'{e}airement disjointes ($K/m$ est modulaire), $B$ est aussi une base de $m(K^{p^{n+1}})$ sur $m$. Or, par construction,   $\te_{{i+2},n+1}={Z_{i+1}}^{p^{-1}}\te_{{i+1},2(n+1)}+{\te_{{i+2},n}}^{p^{-1}}$, donc ${\te_{{i+2},n+1}}^{p^{n+1}}={Z_{i+1}}^{p^n}{\te_{{i+1},2(n+1)}}^{p^{n+1}}+{\te_{{i+2},n}}^{p^n}$, avec ${\te_{{i+2},n}}^{p^n}={Z_{i+1}}^{p^{n-1}}{\te_{i+1,2n}}^{p^n}+\cdots + Z_{i+1}{\te_{i+1,2}}^{p}+Z_{i+2}\in m$.  Par identification, ${Z_{i+1}}^{p^n}\in m\cap K^{p^{n+1}}\subseteq K^{p^{n+1}}$, et donc ${Z_{i+1}}^{p^{-1}}\in K$, absurde. D'o\`{u} pour tout $n\in \N^*$, $\te_{{i+1},n}\in m$, ou encore $K_{i+1}\subseteq m$.  D'o\`{u} $m=K$. \cqfd
%%%%%%%%%%%%%%%%%%%%%%%%%%%%%%%%%%%%
%%%%%%%%%%%%%%%%%%%%%%%%%%%%%%%%%%%%%%%%%%%%%%%%%%%%%%%%%%%%%%%%%%%%%%%µµµµµµµ
%%%%%%%%%%%%%%%%%%%%%%%%%%%%%%%%%%%%%%%%%%%%%%%%%%%%%%%%%%%%%%%%%%%%%
%%%%%%%%%%%%%%%%%%%%%%%%%%%%%%%%%%%%%%%%%%%%%%%%%%%%%%%%%%%%%%%%%%%%%%%%%%%%%%%%%%%??????????????????????????????????????????
\section{Extensions $+\infty$-$w_0$-g\'en\'er\'ees}
 %%%%%%%%%%%%%%%%%%%%%
\subsection{$i$-suite}
\begin{df}
Une suite  $k=K_0 \subseteq K_1\subseteq  \dots \subseteq   K_n\subseteq  \dots \subseteq K$ de sous-extensions d'une extension purement ins\'{e}parable $K/k$ est dite $i$-suite dans $K$ si pour tout indice $i$, on a $K_{i+1}/K_i$ est d'exposant non born\'{e}.
\end{df}

Lorsque $K/k$ est d'exposant non born\'{e}, il est imm\'{e}diat que   $k\subseteq  K$ est une $i$-suite  dite triviale, et donc l'ensemble des $i$-suites de $K/k$ est non vide.  Cependant,  $K/k$ n'admet pas de $i$-suite si $K/k$ est d'exposant fini, donc pour \'{e}carter ce cas, on suppose tout au long de cette paragraphe que $K/k$ est d'exposant non born\'{e}. Si de plus, $K/k$ est $q$-finie, on v\'{e}rifie aussit\^{o}t que
$k=K_0 \subseteq  K_1\subseteq  \dots \subseteq   K_n\subseteq  \dots \subseteq K$ est une $i$-suite  si et seulement si il en est de m\^{e}me de $L=L(K_0) \subseteq  L(K_1)\subseteq  \dots \subseteq   L(K_n)\subseteq  \dots \subseteq K$ pour toute sous-extension finie $L$ de $K/k$. En particulier,  $k=K_0 \subseteq  K_1\subseteq  \dots \subseteq   K_n\subseteq  \dots \subseteq K$ est une $i$-suite si et seulement si $k=K_0 \subseteq  rp(K_1/k)\subseteq  \dots \subseteq   rp(K_n/k)\subseteq  \dots \subseteq K$ l'est aussi.
\begin{pro} {\label{pr41}}Toute suite d\'{e}croissante d'une extension $q$-finie est stationnaire.
\end{pro}
\pre  Soient $(K_n/k)_{n\in \N}$ une suite d\'{e}croissante de sous-extensions de $K/k$ et $(F_i/k)_{i\in \N}$ la suite associ\'{e}e \`{a} leurs cl\^{o}tures relativement parfaites. Compte tenu du th\'eor\`eme \textcolor{blue}{\ref{thm4}} et de la proposition \textcolor{blue}{\ref{pr9}}, la suite des entiers $(di(F_i/k))_{i\in \N}$ est  d\'{e}croissante, donc stationnaire \`a partir d'un entier $n_0$, ou encore pour tout $n\geq n_0$, $F_i=F_{n_0}$. En vertu de la monotonie, pour tout $n\geq n_0$,  $[K_{n+1} :F_{n_0}]\leq [K_n :F_{n_0}]$. Autrement dit, la suite des entiers $([K_{n} :F_{n_0}])_{n\geq n_0}$ est d\'ecroissante, donc stationnaire \`a partir d'un entier $e$, ou encore
 pour tout $n\geq e$, $[K_{n} :F_{n_0}]=[K_{e} :F_{n_0}]$. Comme pour tout $n\geq e$, $K_n\subseteq K_e$, on en d\'{e}duit que $K_n=K_e$, pour tout $n\geq e$. \cqfd
\begin{cor} {\label{pr21}} Dans une extension $q$-finie, toute $i$-suite  est  finie.
\end{cor}
\pre Imm\'ediat. \cqfd
\sk

Soit $K/k$ une extension $q$-finie,  on dit que $K/k$ admet une $i$-suite de longueur $n$ si $K$ peut se d\'{e}composer sous-forme d'extensions :  $k=K_0 \subseteq  K_1\subseteq  \dots \subseteq   K_n=K$  telles que $K_{i+1}/K_i$ est d'exposant non born\'{e} pour tout $i\in \{0,\dots,n-1\}$.  D'apr\`{e}s la proposition pr\'{e}c\'{e}dente toute extension $q$-finie d'exposant non born\'{e} admet  une $i$-suite de longueur maximale $n$, elle sera qualifi\'ee de $i$-suite maximale.  En outre, toute $i$-suite peut se prolonger en une $i$-suite de longueur maximale.  Par ailleurs, une $i$-suite de longueur maximale pr\'{e}sente en quelque sorte une certaine forme d'irr\'{e}ductibilit\'{e} dans la mesure o\`{u} entre deux termes cons\'{e}cutifs n'existe aucune extension propre d'exposant non born\'{e}, et donc impossible de d\'ecomposer deux termes cons\'{e}cutifs en $i$-suite de longueur $2$. Il est \`{a} signaler que cette forme d'irr\'{e}ductibilit\'{e}  sera \'{e}tudi\'{e}e avec pr\'{e}cision dans les sections qui suivent.
\begin{rem} En g\'{e}n\'{e}rale, les termes d'une $i$-suite maximale ne sont pas uniques. Toutefois, on peut chercher d'autres formes d'unicit\'{e}, par exemple on peut se demander si une $i$-suite maximale conserve la taille et les exposants des termes \`{a} une permutation pr\`{e}s.
\end{rem}
%%%%%%%%%%%%%%%%%%%%%%%%%%%%%%%%%%%%%%%%%%%
\subsection{Extensions $w_0$-g\'en\'er\'ees}
\begin{df} Une extension purement ins\'eparable est dite $w_0$-g\'en\'er\'ee, s'elle n'admet pas de sous-extensions propres d'exposant non born\'e.
\end{df}

 En d'autres termes, $K/k$ est $w_0$-g\'en\'er\'ee si toutes les sous-extensions propres  de $K/k$ ont un exposant born\'e. En particulier, si  $K/k$ est $q$-finie, alors $K/k$ est $w_0$-g\'en\'er\'ee si pour toute sous-extension propre $L/k$ de $K/k$, on a $L/k$ est finie, et par suite on retrouve la d\'efinition du J.K Devney cf. \cite{Dev2}. Par ailleurs, la $w_0$-g\'en\'eratrice exprime une certaine forme d'irr\'eductibilit\'e dans la mesure o\`u  $K/k$ est ind\'ecomposable sous  forme d'extensions d'exposant non born\'e. Si de plus $K/k$ est d'exposant non born\'e, on v\'erifie aussit\^ot que :
 \sk

  \begin{itemize}{\it
 \item Toute extension $w_0$-g\'en\'er\'ee  est relativement parfaite.
 \item Toute extension $w_0$-g\'en\'er\'ee et $q$-finie  est modulaire sur une extension finie de $k$.
 \item $K/k$ est $w_0$-g\'en\'er\'ee si et seulement si $k \longrightarrow K$ est une $i$-suite de longueur maximale, et $K/k$ est relativement parfaite.
 \item  Pour toute sous-extension $L/k$ d'exposant born\'e de $K/k$, on a $L(K)/L$ est $w_0$-g\'en\'er\'ee si $K/k$ l'est.}
 \end{itemize}
 \sk

 Le r\'esultat ci-dessous assure l'existence des extensions $w_0$-g\'en\'er\'ees. Plus pr\'ecis\'ement, on a :
\begin{thm} {\label{thm26}} L'ensemble $H$ des sous-extensions d'exposant non born\'e d'une extension $q$-finie $K/k$ d'exposant non born\'e est inductif pour la relation d'ordre $K_1\leq K_2$ si et seulement si $K_2\subseteq K_1$. En particulier, $K/k$ admet une sous-extension $w_0$-g\'en\'er\'ee d'exposant non born\'e.
 \end{thm}
\pre D\'ecoule imm\'ediatement de la proposition \textcolor{blue}{\ref{pr41}}. \cqfd
\begin{pro} {\label{pr42}} Toute extension $q$-finie est compos\'{e}e finie d'extensions $w_0$-g\'en\'er\'ees.
\end{pro}
\pre Le r\'{e}sultat est \'evidemment trivial si $K/k$ est finie. Sinon, d'apr\`{e}s le corollaire \textcolor{blue}{\ref{pr21}}, $K/k$ admet une $i$-suite $k=K_0 \subseteq  K_1\subseteq  \dots \subseteq   K_n=K$  de longueur maximale $n$. N\'{e}cessairement,  $K_i \subseteq   K_{i+1}$ est une $i$-suite de longueur maximale $1$, sinon $K/k$ admet une $i$-suite de la longueur  d\'{e}passant $n$, contradiction.  Par suite, on est amen\'{e} \`{a} d\'{e}montrer le r\'{e}sultat pour $k \subseteq  K$ de longueur maximale $1$. En particulier,   $rp(K/k)$ est irr\'{e}ductible dans la mesure o\`{u} $rp(K/k)/k$ n'admet aucune sous-extension proppre  d'exposant non born\'{e}. Toutefois, d'apr\`{e}s la proposition \textcolor{blue}{\ref{pr9}}, $K/rp(K/k)$ est finie, et par suite $K/k$ est compos\'{e}e finie d'extensions $w_0$-g\'en\'er\'ees. \cqfd
\sk

Dans le cas des extensions modulaires, le r\'esultat suivant montre que la $w_0$-g\'en\'eratrice devient une propri\'et\'e intrins\`eque exclusivement li\'ee aux extensions $q$-finies.
\begin{thm} {\label{thm28}} Pour qu'une extension $w_0$-g\'en\'er\'ee $K/k$ soit $q$-finie il faut et il suffit que $lm(K/k)\not= K$.
\end{thm}

La d\'emonstration de ce th\'eor\`eme fait appel au r\'esultat suivant :
\begin{lem} {\label{tlem28}} Soit $K/k$ une extension purement ins\'{e}parable d'exposant non bor\-n\'{e} et de degr\'e d'irratinalit\'e infini. Si $K/k$ est relativement parfaite et modulaire, alors $K/k$ contient une sous-extension propre $L/k$ d'exposant non born\'{e} et modulaire.
\end{lem}
\pre On va construire par r\'{e}currence une suite strictement croissante $(K_n/k)_{n\geq 1}$ de sous-extensions modulaires d'exposant $n$ de $K/k$. Comme $K/k$ est relativement parfaite, d'apr\`{e}s la proposition \textcolor{blue}{\ref{thm13}} et le corollaire \textcolor{blue}{\ref{acor1}}, pour tout $n\geq 1$, $ di(k^{p^{-n}}\cap K/k)=di(k^{p^{-1}}\cap K/k)=di(K/k)$ et $k^{p^{-n}}\cap K/k$ est \'equiexponentielle d'exposant $n$. Soit $G_1$ une $r$-base de $k^{p^{-1}}\cap K/k$, il en r\'{e}sulte que  $k^{p^{-1}}\cap K\simeq \otimes_k (k(a))_{a\in G_1}$. Choisissons un \'{e}l\'{e}ment $x$ de $G_1$, comme $|G_1|$ est infini,  il existe un sous-ensemble fini $G'_1$ de $G_1$ tel que $x \not \in k(G'_1)$. Posons $K_1=k(G'_1)$, il est clair que $K_1/k$ est modulaire. Supposons qu'on a construit une suite de sous-extensions finies  $k\subseteq K_1\subseteq K_2\subseteq \dots K_n$ de $K/k$ telle que
\sk

 \begin{itemize}{\it
\item[\rm{(1)}] Pour tout $i\in \{1,\dots,n\}$, $K_i/k$ est modulaire.
\item[\rm{(2)}] Pour tout $i\in \{1,\dots,n\}$, $o_1(K_i/k)=i$.
\item[\rm{(3)}] $x\not\in K_n$.
 }
 \end{itemize}
\sk

Soit $G_{n+1}$ une $r$-base de $k^{p^{-n-1}}\cap K/k$, d'apr\`{e}s la proposition \textcolor{blue}{\ref{pr26}}, $k^{p^{-n-1}}\cap K\simeq \otimes_k (k(a))_{a\in G_{n+1}}$. Comme $o_1(K_{n}/k)=n$, on en d\'{e}duit que $K_n\subseteq k^{p^{-n-1}}\cap K$. Or $K_n/k$ est finie et $|G_{n+1}|$ est infini, donc il existe une partie finie $G'_{n+1}$ de $G_{n+1}$ telle que $K_n\subseteq k(G'_{n+1})$. Deux cas peuvent se produire :
\sk

1-ier cas si $x\not\in k(G'_{n+1})$, alors $K_{n+1}=k(G'_{n+1})$ convient.
\sk

2-i\`{e}me cas si $x\in k(G'_{n+1})$, comme $k^{p^{-n-1}}\cap K\simeq \otimes_k (\otimes_k (k(a))_{a\in G'_{n+1}})\otimes_k (\otimes_k (k(a))_{a\in G_{n+1}\setminus G'_{n+1}})$, donc $x\not\in k(G_{n+1}\setminus G'_{n+1})$ ; sinon puisque $k(G'_{n+1})$ et $k(G_{n+1}\setminus G'_{n+1})$ sont $k$-lin\'{e}airement disjoints, alors $x \in k(G'_{n+1})\cap k(G_{n+1}\setminus G'_{n+1})=k$, absurde.  Soit $y$ un \'{e}l\'{e}ment  de $G_{n+1}\setminus G'_{n+1}$, ($y$ existe car $|G_{n+1}|$ est infini et $|G'_{n+1}|$ est fini).  Notons $K_{n+1}=K_n(y)$, on v\'{e}rifie aussit\^ot que :
\sk

 \begin{itemize}{\it
\item $K_{n+1}/k$ est finie, et $o_1(K_{n+1}/k)=o(y,k)=n+1$.
\item $K_{n+1}\simeq K_n\otimes_k k(y)$, (application de la transitivit\'{e} de la lin\'{e}arit\'{e} disjointe de $k(G'_{n+1})$ et $k(G_{n+1}\setminus G'_{n+1})$), et comme $K_n/k$ est modulaire, d'apr\`{e}s la proposition \textcolor{blue}{\ref{apr1}}, $K_{n+1}/k$ est modulaire.
\item $x \not\in K_{n+1}$, sinon comme $k^{p^{{-n-1}}}\cap K\simeq k(G'_{n+1})\otimes_k k(G_{n+1}\setminus G'_{n+1})\simeq K_n(G'_{n+1})\otimes_{K_n} K_n(G_{n+1}\setminus G'_{n+1})$, alors $x\in k(G'_{n+1})\cap K_n(y)\subseteq K_n(G'_{n+1})$ $\cap K_n (G_{n+1}\setminus G'_{n+1})=K_n$, absurde.
 }
 \end{itemize}
\sk

\noindent D'o\`u $K_{n+1}/k$ convient, et par suite $L=\di \bigcup_{i\geq 1}K_i$ satisfait les conditions du lemme ci-dessus.\cqfd
\sk

\noindent{}{\bf Preuve du th\'eor\`eme \textcolor{blue}{\ref{thm28}}}. La condition n\'{e}cessaire r\'{e}sulte imm\'{e}diatement du th\'{e}or\`{e}me \textcolor{blue}{\ref{thm27}}. Inversement, soit  $m$ la plus petite sous-extension de $K/k$ telle que $K/m$ est modulaire. Comme $K/k$ est $w_0$-g\'en\'er\'ee et $m\not=K$, alors $m/k$ admet un exposant fini que l'on note $e$, et d'apr\`{e}s le lemme \textcolor{blue}{\ref{tlem28}} ci-dessus, $K/m$ sera $q$-finie. Dans la suite, pour tout $n\in\N^*$, posons $K_n=m^{p^{-e-n}}\cap K$ et $di(K/m)=l$. Soit $G_n$ une $r$-base de $K_n/m$, compte tenu de la proposition \textcolor{blue}{\ref{thm13}} et le corollaire \textcolor{blue}{\ref{acor1}}, $|G_n|=l$ et $o_1(K_n/m)=e+n$. Par ailleurs, on a $k({K_n}^{p^e})=k(m^{p^e}, {G_n}^{p^e})=k({G_n}^{p^{e}})$. En outre, $di(k({K_n}^{p^e})/k)\leq l$, et $o_1(k({K_n}^{p^e})/k)\geq o_1(m({K_n}^{p^e})/m)=n$. En particulier, l'extension $L=\bigcup k({K_n}^{p^e})$ est d'exposant non born\'{e}, mais comme $K/k$ est $w_0$-g\'en\'er\'ee, on obtient $K=L$. Toutefois, en vertu de la proposition \textcolor{blue}{\ref{pr6}}, $di(L/k)=\di\sup_{n\in\N}(di(K_n/k))\leq l$, donc $K/k$ est $q$-finie.
%\cqfd
%%%%%%%%%%%%%%%%%%%%%%%%%%%%%%%%%%%%%%%%
\begin{cor} {\label{tcor31}} Toute extension modulaire et $w_0$-g\'en\'er\'ee est $q$-finie.
\end{cor}

Compte tenu du th\'eor\`eme ci-dessus et dans le but d'\'etendre la notion de $w_0$-g\'en\'eratrice, on adopte le point de vue suivante :
\subsection{G\'en\'eralisation d'une extension $w_0$-g\'en\'er\'ee}
\begin{df} Soit $j$ un entier naturel non nul. Une extension purement ins\'eparable $K/k$ est dite $j$-$w_0$-g\'en\'er\'ee si  $K/k$ n'admet pas de sous-extensions propres d'exposant non born\'e et degr\'e d'irrationalit\'e  inf\'erieur  \`a $j$.
\end{df}

Autrement dit, toute extension propre de $K/k$ dont le degr\'e d'irrationalit\'e ne d\'epasse pas $j$ strictement est d'exposant fini.
\begin{df}Une extension purement ins\'eparable $K/k$ est dite $+\infty$-$w_0$-g\'e\-n\-\'e\-r\'ee si  pour tout $j\in \N^*$, $K/k$ est $j$-$w_0$-g\'en\'er\'ee.
\end{df}

Il est clair que
toute extension $w_0$-g\'en\'er\'ee  est $+\infty$-$w_0$-g\'en\'er\'ee. Notamment, ces deux notions co\^{\i}ncident dans le cas de la $q$-finitude.
Toutefois, pour \'eviter la non-contradiction, la construction d'un exemple d'extension $+\infty$-$w_0$-g\'en\'er\'ee de degr\'e d'irrationalit\'e infini n\'ecessite les r\'esultats suivants :
\begin{thm} {\label{thm29}} Etant donn\'ee  une extension purement ins\'{e}parable $K/k$ relativement parfaite et modulaire ; et soit $L/k$ une sous-extension propre finie de $K/k$. Si $K/L$ est modulaire, alors pour tout entier $n>e=o_1(L/k)$, $k^{p^{-n}}\cap K/k(L^{p^{e-1}})$ est modulaire. En particulier, $K/k(L^{p^{e-1}})$ est modulaire.
\end{thm}

Pour la preuve de ce th\'{e}or\`{e}me, on se servira des r\'{e}sultats suivants.  D'abord pour tout $n\in \N$, on pose $K_n=k^{p^{-e-n}}\cap K$ et $L_n=L^{p^{-n}}\cap K$.
\begin{lem} {\label{lem18}} Sous les m\^{e}mes hypoth\`{e}ses du th\'{e}or\`{e}me ci-dessus. Pour tout $n\in \N^*$, il existe deux sous-extensions $N$ et $M$ de $K_n/k$ v\'{e}rifiant :
\sk

 \begin{itemize}{\it
\item $L\subseteq k(N^{p^n})$, avec $N/k$ est finie.
\item $K_n\simeq M\otimes_k N\simeq (M\otimes_k L)\otimes_L N$. En outre, $M/k$ et $N/k$ sont \'equiexponentielles d'exposant $n+e$
\item $L(M)/ L(M^p)$ et $L({L_{n+e}}^{p})/L(M^p)$ sont $L(M^p)$-lin\'{e}airement disjointes.
\item $L_{n+e}/L(M)$ est modulaire avec $di(L_{n+e}/L(M))=di(K_n/M)=di(N/k)$.
 }
 \end{itemize}
\end{lem}
\pre Puisque $L/k$ est d'exposant $e$, donc $L\subseteq k^{p^{-e}}\cap K$. D'o\`{u} $L\rightarrow L^{p^{-n}}\cap K\rightarrow K_n \rightarrow L_{n+e}$. Soit $G$ une $r$-base de $K_n/k$, comme $K/k$ est relativement parfaite et modulaire, alors d'apr\`{e}s la proposition \textcolor{blue}{\ref{thm13}}, $K_n/k$ est \'equiexponentielle d'exposant $e+n$. En outre, $K_n\simeq \otimes_k(k(a))_{a\in G}$, et donc  $K_0= k({K_n}^{p^n})\simeq \otimes_k (k(a^{p^n}))_{a\in G}$. Or, $L/k$ est finie, et $L\subseteq K_0$, donc il existe une partie finie $G_1$ de $G$ telle que $L\subseteq k({G_1}^{p^n})$. Notons le compl\'{e}mentaire de $G_1$ dans $G$ par $G_2$, ($G_2=G\setminus G_1$), et d\'esignons respectivement par $M$ et $N$ les corps $k(G_2)$ et $k(G_1)$. On v\'erifie aussit\^ot que :
\sk

 \begin{itemize}{\it
\item $K_n\simeq M\otimes_k N\simeq (M\otimes_k L)\otimes_L N $.
\item $M$ et $N$ sont \'equiexponentielle d'exposant $e+n$.
 }
 \end{itemize}
\sk

\noindent En particulier, pour tout $x\in G_2$, $o(x,L(G_2\setminus \{x\}))=n+e$ ; et par suite s'il existe $x\in G_2$ tel que $x\in L({L_{n+e}}^p)(G_2\setminus \{x\})$, on aura $n+e=o(x,L(G_2\setminus \{x\}))\leq o_1(L({L_{n+e}}^p)(G_2\setminus \{x\})/L(G_2\setminus \{x\}))\leq o_1(L({L_{n+e}}^p)/L)=n+e-1$, c'est une contradiction. D'o\`u, $G_2$ est $r$-libre sur $L({L_{n+e}}^p)$, ou encore $L(M)/L(M^p)$ et $L({L_{n+e}}^p)$ sont $L(M^p)$-lin\'{e}airement disjointes. D'apr\`es le th\'eor\`eme de la $r$-base incompl\`ete, il existe une partie $G_3$ de  $L_{n+e}$ telle que $G_2\cup G_3$ est une $r$-base de $L_{n+e}/L({L_{n+e}}^p)$, compte tenu de la proposition \textcolor{blue}{\ref{pr5}},  $G_2\cup G_3$ est aussi un $r$-g\'en\'erateur minimal de $L_{n+e}/L$. Puisque $K/L$ est modulaire et relativement parfaite, donc $L_{n+e}\simeq \otimes_L (L(a))_{a\in G_2\cup G_3}\simeq (L\otimes_k M)\otimes_L (\otimes_L (L(a))_{a\in G_3})\simeq M\otimes_k (\otimes_L (L(a))_{a\in G_3})$. D'o\`{u} $L_n(M)\simeq M\otimes_k (\otimes_L (L(a^{p^e}))_{a\in G_3}\simeq (M\otimes_k L)\otimes_L (\otimes_L (L(a^{p^e}))_{a\in G_3}\subseteq K_n$ et $K_n\simeq  M\otimes_k N\simeq (M\otimes_k L)\otimes_L N\subseteq L_{n+e}$. D'une part,  comme $N/k$ est \'equiexponnetielle d'exposant $n+e$ et $L\subseteq k(N^{p^n})$, on aura $|G_1|=di(N/k)=di(N/k(N^{p^n}))\leq di(N/L)\leq di(N/k)$, et donc $di(N/L)=|G_1|$. D'autre part, en vertu du th\'eor\`eme \ref{thm6} et du corollaire \ref{cor4}, on a
$|G_3|=di(L_n(M)/L(M))\leq di(K_n/L(M))=di(N/L)$ et $di(K_n/L(M)) \leq di(L_{n+e}/$ $L($ $M))=|G_3|$, (car $K_n\subseteq L_{n+e}$). Par suite, on aura $|G_3|=|G_1|=di(N/k)$.\cqfd
\sk

Comme $K_n\simeq M\otimes_k N\simeq (M\otimes_{k} L)\otimes_L N$ et $K_n/k$ \'equiexponentielle d'exposant $n+e$, on v\'erifie aussit\^ot que :
\sk

 \begin{itemize}{\it
\item Pour tout $i\in\{1,\dots,n\}$, $k({K_n}^{p^i})=K_{n-i}=k(M^{p^i})\otimes_k k(N^{p^i})$.
\item Pour tout $i\in\{1,\dots,n\}$, $M({K_n}^{p^i})=M(K_{n-i})=M \otimes_k k(N^{p^i})$. En particulier, pour tout $i\in\{1,\dots,n\}$, $M(K_i)/M$ est \'equiexponentielle d'exposant $e+i$ et $di(M(K_i)/M)=di(N/k)$.
\item $L_{n+e}/L(M)$ est \'equiexponentielle d'exposant $n+e$. En outre, $L_{n+e}/L(M)$ est modulaire.
 }
 \end{itemize}
\sk

Dans la suite, on pose $di(N/k)=j$, et d\'esignons par $s$ le plus grand entier tel que $o_s(L/k)=o_1(L/k)=e$.
\begin{lem} {\label{lem19}} Sous les conditions ci-dessus, pour tout $n\in \N^*$, on a :
\sk

 \begin{itemize}{\it
\item[\rm(i)] Pour tout $i\in \{0,\dots, n-1\}$, $di(M({K_n}^{p^i})/L(M))=di(N/k)$.
\item[\rm(ii)] $di(M({K_n}^{p^n})/L(M))=di(M(K_0)/L(M))=j-s$.
 }
 \end{itemize}

\noindent
En particulier,  pour tout $r\in\{j-s+1,\dots,j\}$, $$o_r(K_n/L(M))=o_{j-s+1}(K_n/L(M))=n.$$
\end{lem}
\pre Soit $\{\al_1,\dots,\al_m\}$ une $r$-base canoniquement ordonn\'{e}e de $L/k$, donc $k\rightarrow k(\al_1,\dots,\al_s)\rightarrow L\rightarrow  K_0\rightarrow K_n$. Soit $B$ une $r$-base  de $M(K_0)/M(L)$, donc $M(K_0)=M(\al_1,\dots, \al_m, B)$. Or, $L(M)\simeq L\otimes_k M$, donc $M(\al_1,\dots,\al_s)/M$ est \'equiexponentielle d'exposant $e$. Compl\'etons le syst\`eme $\{\al_1,\dots,\al_s\}$ en une $r$-base  de $M(K_0)/M$ par une partie $C$ de $K_0$. En particulier, on aura $|B|=di(M(K_0)/L(M))\leq di(M(K_0)/M(\al_1,\dots, \al_s))$ $= |C|=j-s$. Par ailleurs, pour tout $r\in \{s+1,\dots, m\}$, $o(\al_r,k(\al_1,\dots,\al_s))<e$,  ainsi par application de l'algorithme de la compl\'etion des $r$-bases $M(K_0)=M(\al_1,\dots,\al_s,B)$ ; et donc $B$ est une $r$-base de $M(K_0)/M(\al_1,\dots,\al_s)$. D'o\`u, $di(M(K_0)/M(L))=j-s=|B|$.
De m\^eme, on a $L(M)({K_n}^{p^{n-1}})=K_1(M)$ et $L(M)({K_n}^{p^n})=M(K_0)$. Comme $K_n\simeq
L(M)\otimes_L N$, et donc $M({K_n}^{p^{n-1}})\simeq L(M)\otimes_L L(N^{p^{n-1}})$, il en r\'{e}sulte que $di(M(L)({K_n}^{p^{n-1}})/M(L))=di(M(K_1)/M(L))=di(L(N^{p^{n-1}})/L)$. Or, $N/k$ est \'equiexponentielle d'exposant $n+e$ et $L\subseteq k(N^{p^{n}})$, en vertu du th\'eor\`eme \textcolor{blue}{\ref{thm4}} et le corollaire \textcolor{blue}{\ref{cor4}}, on aura $j=di(k(N^{p^{n-1}})/k(N^{p^{n}}))\leq di(L(N^{p^{n-1}})/L)\leq di(k(N^{p^{n-1}})/k)=j$. Par suite, $di(M(K_1)/M(L))$ $=j$. D'o\`u, d'apr\`{e}s le lemme \textcolor{blue}{\ref{lem2}},  pour tout $r\in\{j-s+1,\dots,j\}$, $o_r(K_n/L(M))=o_{j-s+1}(K_n/L(M))=n$.%\cqfd
\sk

\noindent{}{\bf Preuve du th\'{e}or\`{e}me \textcolor{blue}{\ref{thm29}}.} Tout au long de cette d\'emonstration,  on se servira des notations pr\'{e}c\'{e}dentes. D'abord, pour tout $n\in \N^*$, on a :
\sk

 \begin{itemize}{\it
\item[\rm(1)] $K_n\subseteq L_{n+e}$.
\item[\rm(2)] $L_{n+e}/L(M)$ est modulaire, avec $di(L_{n+e}/L(M))=j=di(K_n/L(M))$.
\item[\rm(3)] $K_n\simeq L(M)\otimes_L N$.
 }
 \end{itemize}

\noindent En vertu de la proposition \textcolor{blue}{\ref{ajpr1}}, il existe une $r$-base canoniquement ordonn\'{e}e $\{a_1,\dots,a_j\}$ de $K_n/L(M)$ telle que $K_n\simeq L(M)\otimes_L L(a_1,\dots, a_{j-s})\otimes_L L(a_{j-s+1})\otimes_L \dots\otimes_L L(a_j)$ ; et donc pour tout $i\in \{j-s+1,\dots, j\}$, ${a_i}^{p^n}\in L$. Soit $\{\al_1,\dots, \al_m\}$ une $r$-base de $L/k$, donc $K_n=M(\al_1,\dots,\al_m,a_1,\dots,a_j)$. Comme $o(\al_i,k)\leq e$ pour tout $i\in \{1,\dots, m\}$ et $K_n/M$ est \'equiexponentielle d'exposant $n+e$, d'apr\`{e}s l'algorithme de la compl\'{e}tion des $r$-bases et le lemme \textcolor{blue}{\ref{lem18}}, $K_n\simeq M(a_1,\dots,a_j)\simeq M\otimes_k k(a_1)\otimes_k\dots\otimes_k k(a_j)$. Or, $K_n/k$ et $K_n/M$ sont \'equiexponentielles  d'exposant $n+e$,  donc $k({a_{j-s+1}}^{p^n},\dots,{a_j}^{p^n})/k$ est \'equiexponentielle d'exposant $e$. D'autre part, $ k({a_{j-s+1}}^{p^n},\dots,{a_j}^{p^n})\subseteq L$,  donc en compl\'{e}tant  ce syst\`eme en
 une $r$-base canoniquement ordonn\'{e}e de $L/k$, on obtient $k(L^{p^{e-1}})=k({a_{j-s+1}}^{p^{n+e-1}},$ $\dots,{a_j}^{p^{n+e-1}})$. Par suite, en vertu de la proposition \textcolor{blue}{\ref{pr24}}, on aura $K_n\simeq M\otimes k(a_1)\otimes_k \dots \otimes_k k(a_j)\simeq M \otimes_{k(L^{p^{e-1}})} k(L^{p^{e-1}})$ $(a_1)\otimes_{k(L^{p^{e-1}})} \dots \otimes_{k(L^{p^{e-1}})} k(L^{p^{e-1}})(a_j)$, avec $M/k$ est modulaire. D'apr\`es la proposition \textcolor{blue}{\ref{apr1}}, on en d\'eduit que $K_n/ k(L^{p^{e-1}})$ est aussi modulaire.\cqfd
 %%%%%%%%%%%%%%%%%%%%%%%%
 \begin{lem} {\label{lem20}} Soit $K/k$ une extension purement ins\'{e}parable \'equiexponentielle d'exposant $n>1$, et telle que $k\not\subseteq K^p$. Soit $L/k$ une sous extension propre de $k^{p^{-1}}\cap K$. Il existe une extension $K'/K$ v\'{e}rifiant les conditions ci-dessous :
 \sk

 \begin{itemize}{\it
\item[\rm(1)] $di(K/k)=di(K'/k)$
\item[\rm(2)] $K'/k$ est \'equiexponentielle d'exposant $n+1$.
\item[\rm(3)] $K'/L$ n'est pas modulaire.
 }
 \end{itemize}
\end{lem}
\pre 1-ier cas : si $K/L$ est non modulaire, alors $K'=K^{p^{-1}}$ convient.
\sk

2-i\`{e}me cas : si $K/L$ est modulaire, d'apr\`{e}s le th\'{e}or\`{e}me \textcolor{blue}{\ref{thm12}}, il existe une $r$-base $G$ de $K/k$ telle que $G_1=\{(a^{p^{o(a,L)}})_{a\in G}|\, o(a,L)<n\}$ est une $r$-base modulaire de $L/k$ et $K\simeq \otimes_L (\otimes_L L(a))_{a\in G}$. Puisque $L\subseteq k^{p^{-1}}\cap K$ et $K/k$ est \'equiexponentielle d'exposant $n$, alors pour tout $a\in G$, on a  $o(a,k^{p^{-1}}\cap K)=n-1\leq o(a,L) \leq o_1(K/k)=n$. Il en r\'esulte que $G_1=\{a\in G$ tel que $o(a,L)=n-1\}$ ; et par suite $K\simeq \otimes_L (L(a))_{a\in G_1}\otimes_L (\otimes_L (L(a))_{a\in G\setminus G_1}$. N\'{e}cessairement,  $G\setminus G_1$ et $G_1$ sont non vides, sinon $k^{p^{-1}}\cap K=L$ ou $L=k$, contradiction avec le fait que $L$ est un sous-corps propre de $k^{p^{-1}}\cap K/k$. Soient $\al \in G\setminus G_1$ et $\be \in G_1$. Comme $k\not\subseteq K^p$, il existe $t\in k$ tel que $t\not\in K^p$. Notons $G'=(a^{p^{-1}})_{a\in G\setminus \{\be\}}\cup \{t^{p^{-1}}{\al}^{p^{-1}}+{\be}^{p^{-1}}\}$ et $K'=k(G')$. On v\'{e}rifie aussit\^{o}t que :
\sk

 \begin{itemize}{\it
\item $K'/k$ est \'equiexponentielle d'exposant $n+1$.
\item $K\subseteq K'$, et $di(K/k)=di(K'/k)$.
 }
 \end{itemize}
\sk

 Si $K'/L$ est modulaire, alors ${K'}^{p^n}$ et $L$ sont $L\cap {K'}^{p^n}$-lin\'{e}airement disjointes. Comme ${\al}^{p^{n-1}}\not\in L$, ou encore $(1,{\al}^{p^{n-1}})$ est $L$-libre, donc $(1,{\al}^{p^{n-1}})$ est en particulier $L\cap {K'}^{p^{n}}$-libre. Compl\'{e}tons ce syst\`{e}me en une base $B$ de ${K'}^{p^n}$ sur $L\cap {K'}^{p^n}$. Compte tenu de la lin\'{e}arit\'{e} disjointe,  $B$ est aussi une base de $L({K'}^{p^n})$ sur $L$. Or, $({\al}^{p^{-1}}t^{p^{-1}}+{\be}^{p^{-1}})^{p^{n}}=t^{p^{n-1}}{\al}^{p^{n-1}}+{\be}^{p^{n-1}}$, par identification on aura $t^{p^{n-1}}\in k\cap {K'}^{p^n}$, et donc $t^{p^{-1}}\in k^{p^{-1}} \cap K' = k^{p^{-1}}\cap K \subseteq K$, c'est une contradiction avec le fait que $t\not\in K^p$. Il en r\'{e}sulte que $K'/L$ est non modulaire. %\cqfd
\begin{lem} {\label{lem21}} Etant donn\'{e}s un corps $k$  de caract\'{e}ristique $p>0$, et $\Om$ une cl\^{o}ture alg\'{e}brique de $k$. Soit $H$ l'ensemble des sous-extensions finies de $\Om/k$.  Si $k$ est d\'{e}nombrable, il en est de m\^{e}me de $\Om$ et $H$.
\end{lem}
\pre   Dans  $\Om$ on d\'{e}finit la relation $\sim$ de la fa\c{c}on suivante : $\al\sim \be$ si et seulement si $irr(\al, k)=irr(\be ,k)$, o\`u $irr(\al, k)$ et $irr(\be ,k)$ sont respectivement les polyn\^omes minimals sur $k$ de $\al$ et $\be$. On v\'erifie imm\'ediatement que  $\sim$ est une relation d'\'{e}quivalence. Soit $E$ un syst\`eme de repr\'esentants dans $\Om$ de cette relation (on peut choisir les \'el\'ements de $E$  parmi les  racines de tous les polyn\^{o}mes irr\'{e}ductibles unitaires de telle mani\`{e}re que chaque polyn\^{o}me sera identifi\'{e} par une et une seule racine, c'est-\`{a}-dire par un \'{e}l\'{e}ment de $E$). D'o\`{u}, $\Om=\di\bigcup_{a\in E} \overline{a}$. Comme les racines d'un polyn\^{o}me sont finies, donc pour tout $a\in E$, $|\overline{a}|$ est fini. De m\^eme, on a $k[X]$ est d\'{e}nombrable, en particulier $E$ l'est aussi ; et par suite $\Om$ est d\'{e}nombrable (cf. \cite{N.B}, III,  p.49, corollaire \textcolor{blue}{3}). Dans la suite, pour tout $n\in \N^*$, notons $H_n=\{L\in H$ tel que $L/k$ est  engendr\'ee par au plus $n$ \'el\'ements$\}$. Il est clair que :
\sk

 \begin{itemize}{\it
\item[$\bullet$] L'application \begin{eqnarray*}
 \Om^n & \longrightarrow & H_n,\\
(\al_1,\dots,\al_n) & \longmapsto  & k(\al_1,\dots,\al_n),
\end{eqnarray*}
 est surjective, donc $|H_n|\leq |\Om^n|$ ; et par suite $H_n$ est d\'{e}nombrable.
\item[$\bullet$] Comme $H=\di\bigcup_{n\geq 1} H_n$, alors $H/k$ est d\'{e}nombrable (cf. \cite{N.B}, III,  p.49, corollaire \textcolor{blue}{3}).\cqfd
 }
 \end{itemize}
\sk

Construisons maintenant, une extension $+\infty$-$w_0$-g\'en\'er\'ee de degr\'e d'irrationalit\'e infini. Pour cela,
consid\'erons un corps commutatif d\'enombrable  $k$ de caract\'eristique $p>0$ et de degr\'e d'imperfection infini, et soit $((X_i)_{i\in \N^*},t)$ une famille $p$-libre sur $k$ . Notons $M_1=k(({X_{i}}^{p^{-1}})_{i\in\N^*})$ et $M_2=k(({X_{i}}^{p^{-2}})_{i\in\N^*})$. D\'esignons par $E$ l'ensemble des sous-extensions propres de $M_1$.  Compte tenu du lemme \textcolor{blue}{\ref{lem21}}, $E$ est d\'enombrable, donc on peut pr\'esenter $E$ sous la forme $E=(L_n)_{n\geq 3}$. Par application du lemme \textcolor{blue}{\ref{lem20}}, on construit une suite d'extensions croissantes $(M_n/k)_{n\geq 3}$ v\'erifiant :
\sk

 \begin{itemize}{\it
\item[\rm(i)] $M_n/L_n$ est non modulaire.
\item[\rm(ii)] $M_n/k$ est \'equiexponentielle d'exposant $n$.
 }
\end{itemize}
\sk

Posons $K=\di\bigcup_{i\in\N^*} M_n$.
\begin{thm} {\label{thm30}} L'extension $K/k$ ci-dessus est modulaire et $+\infty$-$w_0$-g\'en\'er\'ee de degr\'e d'irrationalit\'e infini.
\end{thm}

Pour la preuve on se servira en plus du r\'esultat suivant :
%%%%%%%%%%%%%%%%%%%%%%%%%%%%%%%%%%%%%%%%%%%%%%
\begin{lem} {\label{lem22}}
Etant donn\'ee une extension purement ins\'{e}parable, relativement
parfaite et mo\-du\-lai\-re $K/k$. Soient  $S/k$ une sous-extension  relativement parfaite de $K/k$ et $L/k$ une
sous-extension  de $S/k$. Si $S/L$ est mo\-du\-lai\-re, alors $K/L$ est modulaire, (et en particulier, $K/S$ est modulaire).
\end{lem}
\pre cf. {\cite{Che-Fli2}}, p.
\textcolor{blue}{155}, lemme \textcolor{blue}{2.6}. \cqfd
\sk

\noindent{} {\bf Preuve du th\'eor\`eme \textcolor{blue}{\ref{thm30}}}  D'abord, par construction $K/k$ est relativement parfaite de degr\'e d'irrationalit\'e infini. Soit ensuite $S/k$ une sous-extension propre de $K/k$ de taille finie $di(S/k)=j$ ; supposons que $S/k$ est d'exposant non born\'{e}. En vertu du th\'{e}or\`{e}me \textcolor{blue}{\ref{thm26}}, $S/k$ admet une sous-extension $w_0$-g\'en\'er\'ee que l'on note $S'$, en particulier $S'/k$ est relativement parfaite. Posons $L'=lm(S'/k)$, donc $L'/k$ est finie. Soit $L={L'}^{p^{-1}}\cap S'$ ; d'apr\`{e}s la proposition \textcolor{blue}{\ref{pr23}},  $S'/L$ est modulaire. Compte tenu du lemme \textcolor{blue}{\ref{lem22}} ci-dessus, $K/L$ est modulaire. Par application du th\'{e}or\`{e}me \textcolor{blue}{\ref{thm29}}, pour $n$ assez grand, $k_n/k(L^{p^{e-1}})$ est modulaire, o\`{u} $e=o_1(L/k)$. De plus, on a $ k(L^{p^{e-1}})/k$ est finie d'exposant $1$, donc il existe $t\in \N$ tel que $L_t= k(L^{p^{e-1}})$. Or, $k_n=M_n$ pour tout $n\geq 3$, donc $M_m/L_t$ est modulaire pour un entier assez grand $m>t$, et par suite $L_t(M_{m}^{m-t})/L_t$ est modulaire. D'o\`u $M_t/L_t$ est modulaire, c'est une contradiction avec la construction des $M_n$.\cqfd
%%%%%%%%%%%%%%%%%%%%%%%%%%%%%%%%%%%%%%%%%%%%%%%%%%%%%%
%%%%%%%%%%%%%%%%%%%%%%%%%%%%%%%%%%%%%%%%%%%%%%%%%%%%%%
%%%%%%%%%%%%%%%%%%%%%%%%%%%%%%%%%%%%%%%%%%%%%%%%%%%%%%%%%%%%%%%%%%%%%%%µµµµµµµ

{\small
{\em Authors' addresses}:
{\em EL Hassane Fliouet}, Regional Center for the Professions of Education and Training, Agadir, Morocco
 e-mail: \texttt{fliouet@yahoo.fr}.
 }
\end{document}